\documentclass[journal,twoside]{IEEEtran}
\usepackage{url}

\usepackage{breakurl}

\usepackage{lineno, hyperref} 
\hypersetup{ colorlinks=gray}

\ifCLASSINFOpdf
\else
\fi

\usepackage[utf8]{inputenc}
\usepackage{textcomp}
\usepackage{graphicx}
\usepackage{amsmath}
\usepackage[version=4]{mhchem}
\usepackage{siunitx}
\usepackage{longtable,tabularx}
\usepackage{textcomp}
\usepackage{soul} 
\usepackage[normalem]{ulem}
\usepackage{mathrsfs}

\usepackage{subfigure}
\usepackage{subeqnarray}
\usepackage{soul} 
\usepackage{mathrsfs}
\usepackage{amssymb}
\usepackage{algorithmic}
\usepackage{textcomp}
\usepackage{xcolor}
\usepackage{multirow}
\usepackage{comment}
\usepackage[none]{hyphenat} 

\usepackage{tikz}
\usetikzlibrary{shapes,arrows}

\newcommand{\ignore}[1]{}

\graphicspath{ {figures/} }

\def\BibTeX{{\rm B\kern-.05em{\sc i\kern-.025em b}\kern-.08em
    T\kern-.1667em\lower.7ex\hbox{E}\kern-.125emX}}

\begin{document}
%
\title{A Stochastic Switched Optimal Control Approach to Formation Mission Design for Commercial Aircraft}
%
%
%

\author{Mar\'{\i}a Cerezo-Maga\~na, 
        Alberto Olivares, 
        Ernesto Staffetti
\thanks{{Universidad Rey Juan Carlos, Department of Telecommunication Engineering, Fuenlabrada, Madrid, Spain (e-mail: maria.cerezo@urjc.es; alberto.olivares@urjc.es; ernesto.staffetti@urjc.es).}}
}

%
%

\markboth{IEEE Transactions on Aerospace and Electronic Systems}%
{Cerezo-Maga\~na \MakeLowercase{\textit{et al.}}: A Stochastic Switched Optimal Control Approach to
Formation Mission Design for Commercial Aircraft}
%



\maketitle


\begin{abstract}
This paper studies the formation mission design problem for commercial aircraft in the presence of uncertainties. Specifically, it considers uncertainties  in the departure times of the aircraft and in the fuel burn savings for the trailing aircraft.
Given several commercial flights, the problem consists in arranging them in formation or solo flights and finding the trajectories that minimize the expected value of the direct operating cost of the flights. 
The formation mission design problem is formulated as an optimal control problem of a stochastic switched dynamical system and solved using nonintrusive generalized polynomial chaos based stochastic collocation.
The stochastic collocation method converts the stochastic switched optimal control problem into an augmented deterministic switched optimal control problem. 
With this approach, a small number of sample points of the random parameters are used to jointly solve particular instances of the switched optimal control problem. The obtained solutions are then expressed as orthogonal polynomial expansions in terms of the random parameters using these sample points. 
This technique allows statistical and global sensitivity analysis of the stochastic solutions to be conducted at a low computational cost. 
The aim of this study is to establish if, in the presence of uncertainties, a formation mission is beneficial with respect to solo flight in terms of the expected value of the direct operating costs. 
\textcolor{black}{Several numerical experiments have been conducted in which uncertainties on the departure times and on the fuel saving during formation flight have been considered. The obtained results demonstrate that benefits can be achieved even in the presence of these uncertainties.}

\end{abstract}


%
\IEEEpeerreviewmaketitle

\section{Introduction}
\label{sect:introduction}
This paper studies the formation mission design problem for commercial aircraft in the presence of uncertainties, 
considering extended formations, in which the longitudinal distance between aircraft is between 10 and 40 wingspans.
This paper is actually a follow-up to \cite{cerezo2021formation}, where the potential contribution of formation flight to mitigating the environmental impact of aviation and increasing the capacity of the Air Traffic Management (\textsf{ATM}) system is discussed and the formation mission design problem is studied in the absence of uncertainties using deterministic switched optimal control techniques. 

\textcolor{black}{
The solution of the formation mission design problem in the presence of uncertainties is stochastic, i.e, its components are random processes, which can be characterized by their mean and standard deviation functions. 
The expected values and the standard deviations of the latitude and longitude of the optimal trajectories with respect to the expected value of the timing and the expected values and standard deviations of the timing of the trajectory as functions of the expected value of the distance are useful for \textsf{ATM} purposes, such as conflict detection and traffic synchronization. The expected values and standard deviations of the fuel consumption and final time of each flight of the formation mission permit the Direct Operating Costs (\textsf{DOC}) to be estimated in such a way that airlines can establish to what extent a formation mission is economically beneficial. 
}

\textcolor{black}{
The relative contributions of each random variable to the variability of  the latitude and longitude of the optimal trajectories of each aircraft as functions of the expected value of the timing and the variability of  the timing of the trajectories of each aircraft as functions of expected value of the distance 
allow both the \textsf{ATM} staff to establish
which sources of uncertainty must be reduced to increase the
predictability of the trajectories and airlines to determine what
sources of uncertainty must be reduced to decrease the \textsf{DOC}.
}

\textcolor{black}{
The proposed methodology is an effective tool for solving the formation mission design problem in the presence of uncertainties and quantifying the effects of uncertainties in the departure times of the aircraft and in the fuel burn savings of the trailing aircraft on the solution of the problem. 
}

\textcolor{black}{
Formation flight can have a great added value in two major current air transport-related concerns, 
the decarbonisation of the aviation sector and an increase in the (\textsf{ATM}) capacity \cite{chicagoconvention}.
}
The key enabling factors for this concept of operation are a formation control system to keep the aircraft flying in formation at the optimal relative position to optimize the fuel burn savings for the trailing aircraft 
and an \textsf{ATM} system capable of synchronizing flights departing from different airports to ensure that the formation mission occurs as planned. 
In extended formations, due to the great distance between the leader and the follower aircraft, instabilities, such as meandering and external factors, may affect the motion of the vortices. It is therefore important to maintain the relative position between the follower aircraft and the leader's wake vortices precisely, as the fuel burn savings are very sensitive to that relative positioning. This can be done using a formation control system capable of continuously locating the wake vortices, maintaning the aircraft in the optimal relative position, and, in this way, optimizing the fuel savings during the formation flight \cite{capraceetal:2019:wvdatfaff}.

\textcolor{black}{
While technical issues related to maintain an efficient and safe formation flight have been solved in the last years,
some concerns have yet to be addressed.
%
In particular, a major change in airworthiness standards, policies, and procedures is required.
In the 40th International Civil Aviation Organization (\textsf{ICAO}) Assembly, formation flight was proposed as a strategic objective, and the need to develop a new operational concept which includes reduced separations between aircraft allowing formation was established \cite{abeyratne2020outcome}. 
%
}

\textcolor{black}{
Additionally, the main aircraft manufacturers are making a serious effort to implement formation flight. Airbus, the major European aircraft manufacturer, is undertaking an ambitious project called Fello'fly \cite{Airbusweb}, which aims to demonstrate the technical, operational, and economical viability of formation flight for long-haul commercial flights. In November 2021, the first flight test was carried out in transoceanic flights\footnote{https://simpleflying.com/airbus-a350s-bird-like-flight/}. 
In this project, Airbus is collaborating with airlines and air navigation service providers in order to tackle the formation flight challenges. The objective of this project is not only to demonstrate the operational feasibility of formation flight, but also identify safety procedures and standards for transatlantic operations, enabling a controlled entry into service by 2025. Boeing has conducted an extended formation flight study, in which flight tests on Boeing's ecoDemonstrator platform have been carried out \cite{flanzer2020advances}.
%
}

It is well known that uncertainty in flight departure times is among the main causes of trajectory uncertainty, which generates inefficiency in the \textsf{ATM} system \cite{rivasandvazquez:2016:u}. 
However, timing  is a crucial factor in formation missions. Indeed, considering the usual cruise speed of most long-haul commercial aircraft, missing the rendezvous location by, for instance, ten minutes means spatially missing the partner aircraft by 150 km. 
Such cases require catch-up maneuvers, which result in a loss of performance compared to the planned formation mission. 
Thus, uncertainties in both departure times and fuel burn savings for the trailing aircraft must be considered in the formulation of the formation mission design problem and their effects on the resulting aircraft trajectories and the \textsf{DOC}
of the formation mission must be quantified. 
Although aircraft trajectories are not only affected by uncertainty in the initial conditions but also by uncertainties in the aircraft performance models, 
operational uncertainty, and weather uncertainty, mainly wind and storms, this paper only considers uncertainties in the departure times of the aircraft and in the fuel burn savings for the trailing aircraft.  
These uncertain parameters are modeled by means of random variables, which are described by probability distribution functions. The aircraft model is assumed to be known with precision as is the relevant wind field. 

Given several commercial flights, the problem consists in establishing how to organize them in formation or solo flights and in finding the trajectories that minimize the expected value of the \textsf{DOC} of the formation mission. 
This paper only considers formations of up to three aircraft,
because the difficulties in synchronizing more than three flights make formations of more than three aircraft operationally impractical. Moreover, increasing the number of aircraft in the formation asymptotically reduces the benefits obtained from the formation.
Since each aircraft can fly solo or in various positions within a formation, 
the mission is modeled as a stochastic switched dynamical system, in which the
flight modes of the aircraft are described by sets of stochastic ordinary differential equations, the
discrete state describes the combination of flight modes of the individual aircraft, and logical constraints 
establish the switching logic among the discrete states of the system. 
In this paper, the formation mission design problem is formulated as an optimal control problem of a stochastic switched dynamical system.

A  deterministic switched dynamical system is a particular type of hybrid system that consists of several subsystems and a switching law that specifies the active subsystem at each time instant. 
%
%
Deterministic switched dynamical systems are described by both a continuous and a discrete dynamics, in which the transitions among discrete states are not established in advance. In particular, in this paper, each aircraft is assumed to have different flight modes, namely solo flight and flight in different positions within a formation, and their combination is represented by the discrete state of the switched dynamical system, which models their joint dynamic behavior. Each flight mode is represented by different dynamical equations, which may include or not include formation flight benefits in terms of fuel burn savings. Additionally, logical constraints in disjunctive form, based on the stream-wise distance between aircraft, model the switching logic among the discrete states of the system.
Stochastic switched dynamical systems inherit all the features of the deterministic switched dynamical systems.
Additionally, both the discrete and the continuous dynamics are affected by uncertainties.
This general definition can encompass a wide range of stochastic phenomena. 
This paper considers those stochastic hybrid systems in which random variables only affect the continuous dynamics. More specifically, random variables represent uncertain parameter values and uncertain initial conditions.
The probability distribution functions of these random variables are assumed to be known.
 In this case, the continuous dynamics is described by a set of stochastic differential equations but the discrete dynamics is deterministic. The adjective stochastic as used in this paper means that the solutions of the differential equations depend on a vector of random variables. 
A systematic account of recent developments regarding deterministic and stochastic hybrid systems is given in the monography \cite{alwanandliu:2018:tohsdas}.

An optimal control problem of the stochastic switched dynamical system described above
is an optimal control problem in which the continuous dynamics of the system is represented by stochastic 
differential equations, the objective functional is a stochastic functional, and the constraints, which are defined by 
means of stochastic functions, must be satisfied almost surely, i.e., with probability 1. The set of possible exceptions in which the constraints are not satisfied may be non-empty but must have probability 0. 
The adjective stochastic for these functions and the objective functional means in this paper that they depend on a vector of random variables. 
This problem is referred to as the Stochastic Switched Optimal Control Problem (\textsf{SSOCP}).

Modeling, control, and optimal control of stochastic switched systems are addressed in \cite{alwanandliu:2018:tohsdas},  \cite{zhu:2019:uoc},
%
%
%
%
%
%
%
%
in which the random factors acting on the continuous dynamic models of the switched system in each discrete state are represented by some idealized processes, such as the Wiener process, and tools such as stochastic calculus have been employed to obtain solutions.

Another approach to model uncertainties in switched dynamical systems is to treat uncertainties as random variables or random processes and recast the original deterministic switched dynamical system as a stochastic switched dynamical system. This type of stochastic systems are different from those represented by classical stochastic differential equations, where the random inputs are idealized processes.
Consider a stochastic optimal control problem in which only random variables are present in its formulation like in the formation mission design problem studied in this paper.
One of the most commonly used methods to solve optimal control problems of stochastic dynamical systems of this type is the Monte Carlo sampling, in which independent realizations of the random variables are generated based on their probability distributions. For each realization, the optimal control problem becomes deterministic. 
The need for large number of solutions for accurate results can lead to an excessive computational cost, especially for optimal control problems that are already computationally intensive in the deterministic settings, as the mission design problem considered in this paper. The Generalized Polynomial Chaos (\textsf{gPC}) expansion can contribute to alleviating this drawback.
A systematic and coherent presentation of numerical strategies for uncertainty quantification and stochastic computing is given in \cite{xiu:2010:nmfscasma}, with a focus on the methods based on the \textsf{gPC} approach.  

In stochastic collocation methods  the stochastic model equations are satisfied at a discrete set of points, called nodes, in the corresponding random space. 
Polynomial approximation theory is used to locate the nodes strategically to increase the numerical accuracy. 
With this approach, a small number of sample points of the random variables are used to solve jointly particular instances of the \textsf{SSOCP}. The obtained solutions are then expressed as orthogonal polynomial expansions in terms of the random variables using these sample points. 
This is a nonintrusive methodology because the model equations are not altered.
This technique allows statistical and sensitivity analysis of the stochastic solutions to be conducted at a low computational cost. 
Thus, the \textsf{gPC} method converts the \textsf{SSOCP} into an augmented deterministic Switched Optimal Control Problem (\textsf{SOCP}), in which particular instances of the \textsf{SSOCP} are solved together as a single deterministic optimal control problem.
In this paper, the resulting \textsf{SOCP} is solved using the method described in  \cite{cerezo2021formation}.
It is important to point out that the \textsf{gPC} method is applicable to solve the \textsf{SSOCP} when the solution depends smoothly on the random variables. This means that the solutions obtained for each combination of sample points of the random variables must give rise to the same discrete solution, i.e., to the same sequence of discrete states of the switched dynamical system that represents the formation mission. 

Aircraft trajectories are specified by a sequence of geographical coordinates of spatial positions and the timing, i.e., a sequence of time instants at which the corresponding point of the trajectory must be reached. Thus, any delay of the aircraft in reaching a spatial position is considered a temporal deviation from the trajectory as much as a spatial deviation at a given time instant.
The uncertain position of an aircraft at a given time can be described by a region of confidence.
Likewise, the uncertain timing of an aircraft trajectory at a given spatial position can be described by an interval of confidence.
In this paper, the expected values of the geographical coordinates  together with the expected value of timing obtained in the solution of the \textsf{SSOCP} are considered as the reference trajectories to be followed by the aircraft of the formation mission. 

As previously mentioned, the solution of the formation mission design problem in the presence of uncertainties is stochastic, i.e., 
its components are random processes which can be characterized by their mean and standard deviation functions.
Given the \textsf{gPC} expansion, these functions can be directly computed from the coefficients of this expansion \cite{xiu:2010:nmfscasma}.
The expected value coincides with the first coefficient, whereas the variance is the sum of the squares of the other coefficients of the \textsf{gPC} expansion.
The statistical information estimated from the stochastic solutions includes the expected values and standard deviations of the latitude and longitude of the optimal trajectories 
and the other state variables with respect to time, together with the expected values and standard deviations of the timing of the trajectory as functions of the distance.
They also include the expected values and the standard deviations of the arrival times and, in the case of formation flight, the expected values and the standard deviations of the latitude, longitude, and time of the rendezvous and splitting locations. The expected values and standard deviations of the fuel consumption of each flight of the formation mission are also estimated. 
This information is combined to estimate the expected value of the \textsf{DOC}.

Additionally, in the numerical experiment that involves two random variables, a global sensitivity analysis of the stochastic solution is also conducted \cite{saltellietal:2008:gsatp}. 
In this paper, the global sensitivity analysis is based on the Sobol' indices, which enable the determination of what proportions of the variance of the solutions of the switched optimal control problem can be attributed to the different random variables of the switched optimal control model.
The purpose of the sensitivity analysis is to identify the random variables that have more influence on the variability of a component of the stochastic solution with the aim of reducing its variability acting on the source.

\bigskip

\subsection{Previous approaches}

In \cite{shoneetal:2021:aosmiatm},
a review of the literature on stochastic modeling with applications to \textsf{ATM} is provided, including literature on stochastic optimal control.

In \cite{huetal:2005:acpitposcwf},
the problem of aircraft conflict prediction is studied for two-aircraft midair encounters.
First, a model is presented for prediction of the aircraft positions along a time horizon, during which each aircraft is following a prescribed flight plan in the presence of additive wind perturbations on its velocity.
Then, a method for estimating the probability of conflict is proposed.
This method is based on a Markov chain approximation of the stochastic processes that models the aircraft flight. 
In \cite{prandiniandhu:2009:aorafshstacp},
the aircraft conflict prediction problem is formulated as a reachability problem in a stochastic hybrid system framework. Specifically, a switching diffusion model is employed to predict the future positions of an aircraft following a given flight plan, and the probability that the aircraft enters an unsafe region of the airspace is estimated 
using a numerical algorithm for reachability computation.


%
%
%
%
%
%

%
%
%
%
%

In \cite{lietal:2014:artounpc},
an approach to aircraft trajectory optimization in the presence of uncertainties based on \textsf{gPC} is presented, 
where only one random variable has been considered, which
represents an uncertain aerodynamic parameter of the dynamic model of the aircraft. This random variable has been modeled as a uniformly distributed random variable.
In \cite{matsuno2015stochastic},
a stochastic optimal control method based on \textsf{gPC} is developed for determining conflict-free aircraft trajectories under wind uncertainty. 
The random processes that represent the components of the wind speed are approximated as a linear combination of deterministic functions multiplied by independent random variables using the Karhunen-Lo\`eve expansion.

In \cite{cerezo2021formation},
a deterministic switched optimal control method is employed to solve the formation mission design problem. An additional analysis has been conducted to determine how departure delays and fuel burn savings affect the formation flight, in terms of routes, rendezvous and splitting locations and times, and flight times.
In \cite{kentandrichards:2014:afteogdocff},
the impact of ground delays on formation flight is studied using stochastic dynamic programming.

\subsection{Contributions of the paper}

This paper proposes a methodology for the solution of the formation mission design problem in the presence of uncertainties in some of the parameters and boundary conditions of the problem. These uncertainties are represented by random variables characterized by probability density functions.
The formation mission design problem in the presence of uncertainties is formulated as an optimal control problem of a stochastic switched system. 
This problem is solved using an approach based on \textsf{gPC}, in which the stochastic switched optimal control problem is transformed into an equivalent deterministic switched optimal control problems in a higher dimensional state space and solved using a numerical method developed by the same authors \cite{cerezo2021formation}. 
The formulation of the problem based on switched optimal control allows accurate dynamic models of the aircraft and meteorological forecast to be included in the problem formulation. 
This technique permits statistical and global sensitivity analysis of the stochastic solution to be conducted at a low computational cost.

\subsection{Organization of the paper}

The paper is organized as follows. 
The model of the switched dynamical system that represents the formation mission is introduced in Sec.~\ref{sect:model_of_the_system}.
The deterministic and the stochastic switched optimal control problems are stated in sections \ref{sect:the_deterimistic_switched_optimal_control_problem} and \ref{sect:the_stochastic_switched_optimal_control_problem}, respectively. 
The generalized polynomial chaos expansion technique and the method to solve the stochastic  switched optimal control problem are described in Sec.\ref{sect:the_generalized_polynomial_chaos_expansion}.
The results of the numerical experiments are reported and analyzed in Sec.~\ref{sect:numerical_results}.
The outcomes of the sensitivity analysis of the solutions are discussed and interpreted in Sec.\ref{sect:sensitivity_analysis}, 
Finally, conclusions are drawn in Sec.~\ref{sect:conclusions}.

\section{Model of the system}
\label{sect:model_of_the_system}

This section outlines the model of the switched dynamical system that represents the formation mission. 
Further details on the aircraft equations of motions, flight envelope, and wind model employed can be found in \cite[Sect.~II]{cerezo2021formation}.

The mission design problem is studied only in the cruise phase of the flight. 
Thus, a simplified two-degrees-of-freedom point variable-mass dynamic model is considered assuming that all aircraft are in the cruise phase. The motion is restricted to the horizontal plane at cruise altitude over a spherical Earth model. A symmetric flight without sideslip is considered and all the aircraft forces are supposed to be in the plane of symmetry of the aircraft. Wind effects are also considered. All aspects associated with the rotational dynamics are neglected.  
The set of kinematic and dynamic differential-algebraic equations (\textsf{DAE}) that describe the motion of each aircraft of the formation mission are 
\begin{eqnarray} \label{eq:kinematics_dinamics_eqs}
\nonumber  \dot{\phi}(t) &   =  &   \dfrac{V(t) \cdot \cos \chi(t) + {V_{W_N}}(t)}{R_E + h},\\  
  	\dot{\lambda}(t) &   = &   \dfrac{V(t) \cdot \sin \chi(t) + {V_{W_E}}(t)}{ \cos \phi(t) \cdot \left(R_E + h \right) }, \\  
\nonumber  \dot{\chi} (t) &   =  &  \dfrac{L(t) \cdot \sin \mu(t)}{V(t) \cdot m(t)},\\ 
\nonumber   \dot{V}(t) &   =  &   \dfrac{T(t) - D(t)  }{m(t)}, 
\end{eqnarray}
where the state vector has five components: the two dimensional position variables, latitude and longitude, denoted by $\phi$ and $\lambda$, respectively, the heading angle $\chi$, the true airspeed $V$, and the mass of the aircraft $m$. In this set of equations, the control vector has three components: the thrust force $T$, the lift coefficient $C_L$, and the bank angle $\mu$. The normalized version of the \textsf{DAE} system (\ref{eq:kinematics_dinamics_eqs}) is used in this paper.
$L = qS C_L $ is the lift force, 
where $q = \frac{1}{2} \rho V^2$ is the dynamic pressure, $\rho$ is the air density, and $S$ is the reference wing surface area. The aerodynamic drag force is $D = qS C_D $, where $C_D$ is the drag coefficient. 
$h$ is the cruise altitude and ${V_{W_E}}$ and ${V_{W_N}}$ are the components of the wind velocity vector in eastward and northward directions, respectively.
The mass flow rate equation is
\begin{equation}
\dot{m}(t) = - T(t) \cdot  \eta(t), 
\label{mass_differential_eq} 
\end{equation}
where $\eta$ is the thrust specific fuel consumption. The Eurocontrol' s Base of Aircraft Data (\textsf{BADA}), version 3.6 \cite{eurocontrol:2013:umftboad}, has been used to determine this. 
Thus, for aircraft $p$, the state vector is  $x_p(t) = \left( \phi_p(t), \lambda_p(t), \chi_p(t), V_p(t), \right.$
$\left. m_p(t)\right) , \forall p \in  \left\lbrace 1, \ldots, N_a\right\rbrace $, and the control vector is $u_p(t) = \left( T_p(t), C_{L_p}(t), \right.$
$\left. \mu_p(t)\right) ,\forall p \in \left\lbrace 1, \ldots, N_a\right\rbrace $, with $N_a$ the number of aircraft involved in the mission design problem.

Flight envelope constraints model aircraft performance limitations. These constraints are algebraic constraints usually expressed as simple bound constraints that involve, for instance, flight altitude, load factor, and airspeed.
\textsf{BADA} 3.6 is used to model the upper and lower bounds of the flight envelope constraints.
These constraints are included in the set of path constraints of optimal control problems.

Wind information from ERA-Interim \cite{Dee:2011ex}, provided by the European Centre for Medium-Range Weather Forecasts (\textsf{ECMWF}) is used to compute the ${V_{W_E}}$ and ${V_{W_N}}$ components of the wind speed in Eqs.~(\ref{eq:kinematics_dinamics_eqs}). Since this information is given as grid data, analytic functions that approximate the grid data must be determined to include wind information in the stochastic switched optimal control problem used to solve the mission design problem. Radial basis functions \cite{buhmann2003radial} are employed for this purpose. Wind data from April 30, 2019, at 12:00 are used.

The \textsf{DAE} system consisting of Eqs.~(\ref{eq:kinematics_dinamics_eqs}) and Eq.~(\ref{mass_differential_eq}) is applicable to any aircraft in the cruise phase of solo flights. In formation flights, it is necessary to include a modification in the mass flow rate equation (\ref{mass_differential_eq}) for the trailing aircraft.  It is assumed that the reduction in the induced drag of the trailing and intermediate aircraft can be modeled as a percentage of the reduction in fuel consumption, $\mathcal{R}_\text{fuel}$. In this way, the fuel saving factor can be directly incorporated into Eq.~(\ref{mass_differential_eq}) to obtain the mass flow rate equation for the trailing aircraft during formation flight
\begin{equation}  \label{mass_differential_eq_with_reduction}
\dot{m}(t) = - [1 - \mathcal{R}_\text{fuel}] \cdot T(t) \cdot  \eta(t).
\end{equation}
In a two-aircraft formation, this reduction is applied to the trailing aircraft when the distance to the leader aircraft is shorter than 20 wingspans. 
In a three-aircraft formation, this reduction is applied to the trailing and intermediate aircraft when the distances between the leader and intermediate aircraft and between the intermediate and the trailing aircraft are shorter than 20 wingspans.
\textcolor{black}{Notice that the distance of 20 wingspans is used for checking
if aircraft achieve benefits from the formation flight and not for providing flight trajectories strictly accounting for this distance.}

\section{The Deterministic Switched Optimal Control Problem}
\label{sect:the_deterimistic_switched_optimal_control_problem}

Since each aircraft can fly solo or in different positions within a formation, the mission is modeled as a switched dynamical system, in which the discrete state describes the combination of flight modes of the individual aircraft. 

In this section, following \cite[Sect.~III]{cerezo2021formation},
the formulation of the \textsf{SOCP} is introduced for a two-switched dynamical system. 
The dynamical model of a two-switched dynamical system is
\begin{eqnarray} \label{eq:31.1}
\nonumber   \dot{x}_S(t) & =  & f_{v_{S}(t)}(t, x_S(t), u_S(t)), \\  
  x_S(t_I)  & =  & x_I \in \mathbb{R}^n, \\  
\nonumber     x_S(t_F)  & =  & x_F \in \mathbb{R}^n, \\  
\nonumber  v_S (t )  & \in &  \{0, 1\}, \;  t_I \leq t \leq t_F, 
\end{eqnarray}
where the continuously differentiable vector fields, 
$f_0, f_1: \mathbb{R} \times \mathbb{R}^n \times \mathbb{R}^m \rightarrow \mathbb{R}^n$, 
describe the dynamics of the two possible modes of the system.
The control input $u_S(t) \in \Omega \subset \mathbb{R}^m$ is constrained 
to belong, at each time instant,  to the bounded and convex set 
$\Omega$. The binary variable $v_S (t ) \ignore{ \in \{0, 1\} }$ is the 
mode selection variable that specifies which of the two possible 
system modes, $f_0$ or $f_1$, is active. Thus, both $u_S(t)$ and $v_S (t)$
can be regarded as control variables. 
The initial time $t_I$, the final time $t_F$, the initial state $x_S(t_I)$,  
and the final state $x_S (t_F)$ are assumed to be
restricted to a boundary set $B$, namely, $(t_I, x_S (t_I), t_F , x_S (t_F )) \in B = T_I \times B_I \times T_F  \times B_F  \subset \mathbb{R}^{2n+2}$.
The objective functional of the \textsf{SOCP} is
\begin{eqnarray}
&&J_S \left(t, x_S(t) , u_S(t) , v_S(t) \right) =\\
&=& g (t_F, x_F) + \int_{t_I}^{t_F}  F_{v_S(t)} (t, x_S (t), u_S (t)) dt,
\label{eq:31.2}
\end{eqnarray}
where  $g$ is the endpoint cost function defined on a neighborhood of $B$ and $F_0$
and $F_1$ are real-valued continuously differentiable functions that
represent the running cost of the system in each mode. In general, other constraints are included in the formulation of the problem, such as the path constraints, which apply over the whole path or at intermediate points and not only at the end points, and the logical constraints, which 
establish the switching logic among the discrete states of the system. 

The \textsf{SOCP} is stated as follows
\begin{equation}
\min_{u_S \in \Omega, v_S \in \{0,1 \} } J_S \left(t, x_S(t) , u_S(t) , v_S(t) \right),
\label{eq_SOCP_definition}
\end{equation}
subject to Eq.~(\ref{eq:31.1}), endpoint constraints $(t_I, x_S (t_I ), t_F , x_S (t_F)) \in B$, path constraints, and logical constraints.

In the formation mission design problem studied in this paper, not all the dynamic equations are subject to switches. More specifically, the equations of motion associated with the state variables ${\phi}$, ${\lambda}$, ${\chi}$, and ${V}$ do not switch, whereas the aircraft' s mass flow rate equation
does when the aircraft joins or leaves a formation as a trailing or an intermediate aircraft. 
When an aircraft is flying solo or as a leader in the formation it obtains no benefits from the formation flight in terms of fuel savings and Eq.~(\ref{mass_differential_eq}) describes its mass flow. On the contrary, when an aircraft is flying as the intermediate or trailing aircraft in the formation it obtains benefits from the formation flight and its mass flow is described by Eq.~(\ref{mass_differential_eq_with_reduction}). Thus, $f_0$ represents Equations (\ref{eq:kinematics_dinamics_eqs}) and (\ref{mass_differential_eq}), whereas $f_1$ represents Equations (\ref{eq:kinematics_dinamics_eqs}) and (\ref{mass_differential_eq_with_reduction}). 
The functions $F_0$ and $F_1$, which appear in the running cost of Eq.~(\ref{eq:31.2}), are not present in the formulation of the problem. 
The objective {functional} $J_S$ of the switched optimal problem 
that models the formation mission design problem is
\begin{equation}
 J_S = \alpha_t  \sum_{p=1}^{N_a}  t_{{flight}_p}  +  \alpha_f \sum_{p=1}^{N_a}  {{m}_{f_p}}, 
\label{eq:obj_function}
\end{equation}  
where $\alpha_t$ and $\alpha_f$ are the time and the fuel consumption weighting parameters, respectively, $t_{{flight}_p},  \forall p  \in \left\lbrace 1, \ldots, N_a\right\rbrace $, is the total flight time for aircraft $p$, 
and ${{m}_{f_p}} , \forall p  \in \left\lbrace 1, \ldots, N_a\right\rbrace $,  is the fuel consumption for aircraft $p$.

The switching logic among the discrete states of the system is defined by the logical constraints.  
They express the fact that not all transitions between discrete states are possible at any time. 
For example, in a three-aircraft formation mission, the transition between the discrete state 
in which the three aircraft are flying solo and that in which they are flying in a three-aircraft formation
is not permitted. In this case, passing through an intermediate discrete state 
in which two aircraft are flying in formation and the third is still flying solo is mandatory.
Any logic expression can be transformed to conjunctive normal form, i.e., it can be expressed as a conjunction (sequence of ANDs) of one or more clauses, where a clause is a disjunction (OR) of literals,
which represent statements or negations of statements. To incorporate the logical constraints 
in the \textsf{SOCP}, they must be converted into a set of equality or inequality constraints.

Therefore, the formation mission design problem is formulated as an optimal control problem of a switched dynamical system with logical constraints and solved using the techniques described in \cite[Sect.~IV]{cerezo2021formation} and \cite[Sect.~V]{cerezo2021formation}, 
in which binary variables are not used. 
The result is a smooth optimal control problem without binary variables, which reduces the computational complexity of finding a solution. 
This problem is solved using a knotting pseudospectral method as described in \cite[Sect.~VI]{cerezo2021formation}.
The solution of the deterministic optimal control problem is denoted by $z(t)$.
It includes the control and state variables and their timing, from which the rendezvous, splitting, and final times
can be read off.

\section{The Stochastic Switched Optimal Control Problem}
\label{sect:the_stochastic_switched_optimal_control_problem}

Let $\theta = (\theta_1, \theta_2, \ldots, \theta_M)$ be the vector of random variables that represent the random parameters of the \textsf{SOCP}. The presence of random variables in the \textsf{SOCP} converts it into an  \textsf{SSOCP}. In the formation mission design problem studied in this paper, not all the elements of the \textsf{SOCP} contain random variables. Specifically, the fuel consumption reduction factor $\mathcal{R}_\text{fuel}$ in the mass flow rate equation (\ref{mass_differential_eq_with_reduction})  of the aircraft flying in formation as a trailing or an intermediate aircraft  is a random variable. In contrast, the equations of motion associated with the state variables ${\phi}$, ${\lambda}$, ${\chi}$, and ${V}$ do not contain random variables. The departure times of the aircraft are also assumed to be random variables of the \textsf{SOCP}. One of these departure times coincides with $t_I$, the initial time of the formation mission.  These random variables are assumed to be independent and characterized by probability density functions.

The solution of the  \textsf{SSOCP} includes both the state and control variables of the system that represents the formation mission. Due to the presence of the vector random variable $\theta$ in the \textsf{SSOCP}, the solution is stochastic, i.e., it depends on the actual realization of the vector random variable. 
In particular, it is a random function of state or time. Therefore, the solution is actually a random process, 
characterized by the mean function and the corresponding 95\% confidence envelope. 
Specifically, the observed values of a state variable at a given time $t$ is a random variable. 
In this case, the random process is a random function of time. Likewise, the time at which a given value of the state variable is reached is a random variable. This means that the timing of the trajectories, which is defined as the time at which a state is reached, is also a random process. In this case, the random process is a random function of the state. In this paper, the timing of the trajectories is represented as a random function of the orthodromic distance from the departure location.  
For the sake of ease of exposition, the stochastic solution of the \textsf{SSOCP} is denoted in the following sections as a random function of time.
All the random variables observed at specific time instants or states are characterized by their expected values and by the corresponding 95\% confidence intervals, including the arrival times and the fuel consumptions of the aircraft of the formation mission.

The objective functional of the \textsf{SSOCP} is  the expected value of the objective functional (\ref{eq:31.2}) of the \textsf{SOCP}. In the formulation of the \textsf{SSOCP}, the dynamic equations, the path constraints, and the boundary conditions formally depend on the vector random variable $\theta$
and must be satisfied almost surely.




\section{The Generalized Polynomial Chaos Expansion}
\label{sect:the_generalized_polynomial_chaos_expansion}

In this section, following \cite{matsuno2015stochastic}, the \textsf{gPC} expansion is introduced and 
the method for determining the stochastic solution of the \textsf{SSOCP} and computing its statistical information 
is described.

Let $z (t, \theta)$ denote the stochastic solution of the  \textsf{SSOCP}, which includes $x(t, \theta)$ and $u(t, \theta)$. 
The P-th order \textsf{gPC} approximation of $z (t, \theta)$ can be expressed as:
\begin{equation}
z_P (t, \theta) = \sum_{m=1}^M C_m(t) \Phi_m (\theta),
\label{eq:multi_dimensional_expansion_trunc}
\end{equation}
where 
$\theta = (\theta_1, \theta_2, \ldots, \theta_N)$ is a vector of independent random variables,
$C_m(t) $ are the coefficients of the expansion, which can be obtained using either a  nonintrusive or an intrusive approach,
and $\Phi_m (\theta)$ are the multivariate orthogonal polynomial basis functions, which are calculated 
from the $l_i$-th order one-dimensional polynomial basis function $\phi^{(l_i)}(\theta_i)$ of the random variable $\theta_i$
by means of the tensor product rule as follows:
\begin{equation}
\Phi_m(\theta) = \prod_{i=1}^N \phi_i^{(l_i)} (\theta_i).
\label{eq:tensor_product_rule}
\end{equation}
Notice that, in Eq. (\ref{eq:tensor_product_rule}), a unique combination of $l_i, i=1,2, \dots, N$, corresponds to each subscript $m$, which satisfy the condition $\sum_{i=1}^N l_i \leq P$, where $P$ is the maximum degree of the multivariate polynomial $\Phi_m$. The number of tensor product basis functions is $M=\binom{N+P}{N}$. 
The orthonormal polynomials in Eq. (\ref{eq:tensor_product_rule}) satisfy the following orthogonality condition
\begin{equation*}
E[ \phi_i^{(j)} (\theta_i) \phi_i^{(k)} (\theta_i)] = \int \phi_i^{(j)} (\theta_i) \phi_i^{(k)} (\theta_i) \rho_i (\theta_i) d \theta_i =\delta_{jk},
\end{equation*}
where $E$ is the expected value operator,  $\rho_i (\theta_i )$ is the probability density function of the random variable $\theta_i$, and $\delta_{jk}$ is the Kronecker delta function. 

To improve convergence, the choice of the orthonormal polynomials should be made on the basis of the probability density function $\rho_i (\theta_i )$. For instance, the Legendre polynomials are the best choice for uniform random variables, whereas Hermite polynomials are the best option for Gaussian random variables. 
In this paper, the coefficients $C_m$ of the expansion (\ref{eq:multi_dimensional_expansion_trunc})
are calculated using a nonintrusive \textsf{gPC} based stochastic collocation method as follows
\begin{equation}
C_m = E[z(t, \theta) \Phi(\theta)] = \int z(t, \theta) \Phi_m (\theta) \rho(\theta) d \theta,
\label{eq:computation_of_coefficients}
\end{equation}
where $\rho(\theta) = \prod_{i=1}^N \rho_i (\theta_i)$ is the joint probability density function of the vector random variable $\theta$.
A Gaussian quadrature can be used to approximate the integral in Eq. (\ref{eq:computation_of_coefficients}). 
The univariate quadrature rule approximates the polynomial $\phi_i^{(l_i)}(\theta_i)$, $i=1,2,\ldots, N$,
using a set of $q$ collocation points $\theta_i^{(j)}$
and calculating the corresponding set of quadrature weights $\alpha^{(j)}$, $j=1,2,\ldots,q$.
Higher precision in the quadrature can be achieved by increasing the number of collocation points $q$.
Using the tensor product rule, the total number of collocation points is $q^N$. 
This number, which could become large when $N$, the number of random variables in $\theta$, 
grows, can be reduced using the sparse grid quadrature based on the Smolyak rule  \cite{xiu:2010:nmfscasma}.

Let  $\theta^{(j)}$  be the collocation points and $\alpha^{(j)}$, $j=1, \ldots,Q$, the corresponding weights. 
Then, the coefficients of the expansion $C_m$ in Eq. (\ref{eq:computation_of_coefficients}) are approximated by:
\begin{equation}
C_m \approx \sum_{j=1}^Q z(t, \theta^{(j)}) \Phi_m (\theta^{(j)}) \alpha^{(j)},
\label{eq:approximation_of_coefficients}
\end{equation}
where $z(t, \theta^{(j)})$ denotes the solution of the augmented \textsf{SOCP} obtained using the $j$-th collocation point $\theta^{(j)}$
and $\Phi_m (\theta^{(j)})$ represents the multivariate orthogonal polynomial basis function evaluated at the  $j$-th collocation point $\theta^{(j)}$. 
Thus, the \textsf{SSOCP} is converted into an augmented deterministic \textsf{SOCP}, 
in which  particular instances of the \textsf{SSOCP}, which correspond to the collocation points of the random variables, are combined into a single optimal control problem. The resulting augmented \textsf{SOCP} is solved by means of the numerical method described in Sec.~\ref{sect:the_deterimistic_switched_optimal_control_problem}.
The expected value and variance of the stochastic solution $z_P(t, \theta)$ are calculated using the coefficients of the expansion as follows:
\begin{eqnarray}
\text{E}[z_P (t, \theta)]
&=&
C_1(t), \label{eq:mean_gpc}\\
\text{VAR}[z_P (t, \theta)] &=&
\sum_{m=2}^M C_m^2(t) \label{eq:variance_gpc}.
\end{eqnarray}

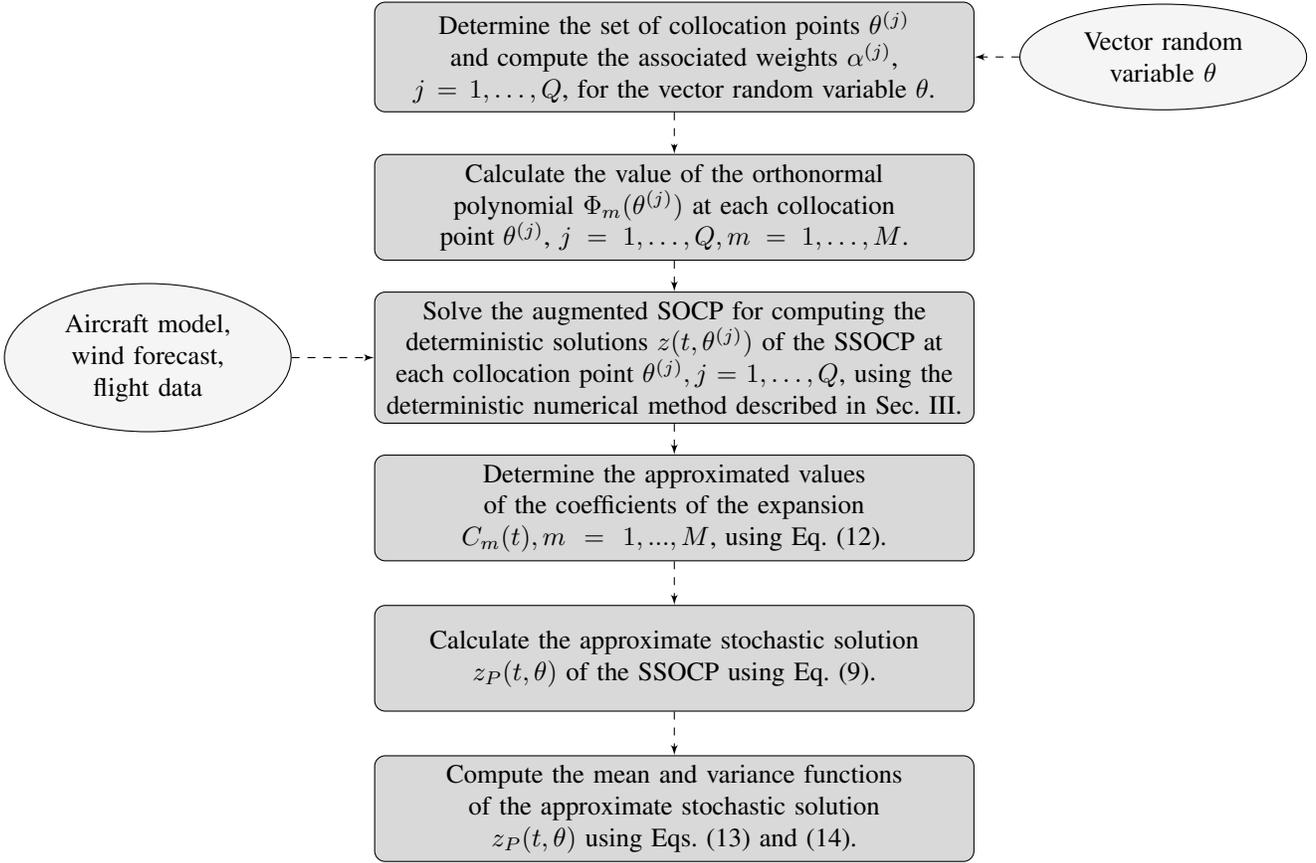
\begin{figure*}
\tikzstyle{decision} = [diamond, draw, fill=gray!30, text width=4.5em, text badly centered, node distance=3cm, inner sep=0pt]
\tikzstyle{block} = [rectangle, draw, fill=gray!30, text width=22em, text centered, rounded corners, minimum height=4em]
\tikzstyle{line} = [draw, -latex'] 
\tikzstyle{cloud} = [draw, ellipse,fill=lightgray!15, node distance=1cm, text width=7em, text badly centered, minimum height=4em]
    
\begin{tikzpicture}[node distance = 2cm, auto]
!\node [block] (uno) {Determine the set of collocation points $\theta^{(j)}$ and compute the associated weights $\alpha^{(j)}$, $j=1, \ldots,Q$, for the vector random variable $\theta$.};
\node [cloud, right of=uno,node distance = 6.5cm] (theta) {Vector random variable $\theta$};
\node [block, below of=uno] (due) {Calculate the value of the orthonormal polynomial $\Phi_m(\theta^{(j)})$ at each collocation point $\theta^{(j)}$, $j=1,\ldots, Q, m=1,\ldots, M$.};
\node [block, below of=due] (tre) {Solve the augmented SOCP for computing the deterministic solutions $z(t,\theta^{(j)})$ of the SSOCP at each collocation point $\theta^{(j)}, j=1, \ldots, Q$, using the deterministic numerical method described in Sec.~III.};
\node [cloud, left of=tre,node distance = 7cm] (parameters) {Aircraft model, wind forecast, flight data};
\node [block, below of=tre] (quattro) {Determine  the approximated values of the coefficients of the expansion $C_m(t), m=1,...,M$, using Eq. (12).};
\node [block, below of=quattro] (cinque) { Calculate the approximate stochastic solution $z_P (t, \theta)$ of the SSOCP using Eq. (9).};
\node [block, below of=cinque] (sei) {Compute the mean and variance functions of the approximate stochastic solution $z_P (t, \theta)$ using Eqs.  (13) and (14).};

\path [line,dashed] (theta) -- (uno);
\path [line,dashed] (parameters) -- (tre);
\path [line,dashed] (uno) -- (due);
\path [line,dashed] (due) -- (tre);
\path [line,dashed] (tre) -- (quattro);
\path [line,dashed] (quattro) -- (cinque);
\path [line,dashed] (cinque) -- (sei);    

\end{tikzpicture}

\caption{
Schematic diagram of the procedure for determining the stochastic solution of the SSOCP and computing its statistical information.
}
\label{fig:schematic_diagram}
\end{figure*}

A complete treatment of the \textsf{gPC} expansion can be found in the monography \cite{xiu:2010:nmfscasma}, whereas the procedure followed to determine the stochastic solution of the \textsf{SSOCP} and to compute its statistical information is schematically represented in Figure~\ref{fig:schematic_diagram}. 
Further details can be found in \cite{lietal:2014:artounpc} and \cite{matsuno2015stochastic}.

\section{Numerical Results}
\label{sect:numerical_results}

In this section,  the results of the following numerical experiments are reported to show the effectiveness of the proposed methodology  to solve the stochastic formation mission design problem: 

\begin{itemize}

\item Experiment A:  Three-aircraft transoceanic mission design with uncertainty in the fuel burn savings.

\item Experiment B: Two-aircraft transoceanic mission design with uncertainty in departure times of the flights.

\end{itemize}

All the experiments involve transoceanic eastbound flights. Wind data
from the ERA-Interim reanalysis database of the \textsf{ECMWF} are used. As mentioned
above, in this paper only the cruise phase is modeled; the rest of
the flight phases are neglected. Thus, the initial and final locations of the cruise 
phase of the flights are assumed to
be the latitudes and longitudes of the departure and arrival airports of
each flight at  cruise altitude. Airbus A330-200 aircraft \textsf{BADA}
models are considered  for each flight. The initial masses of the
aircraft are assigned, as are the initial and final velocities, which
are set at typical cruise values for the selected aircraft models.
The initial heading angles are set to the initial heading angles of the
orthodromic paths between the initial and final locations of each flight.
The time and fuel burn weighting parameters,
$\alpha_t$ and $\alpha_f$, in the \textcolor{black}{objective functional
(\ref{eq:obj_function})} are set to 0.3 and 0.7, respectively, \textcolor{black}{based on 
\cite{cook2004evaluating}, an Eurocontrol study on the delays-related costs}. 

The numerical experiments have been conducted on a 3.6 GHz Intel Core i9 computer with 32 GB RAM. The computational times reported in this section include both the time required to generate the warm-start solution and the time to find the optimal solution of the problem. The warm-start solution is a feasible solution of the problem generated by the NLP solver from an initial guess of the solution.

As previously mentioned, the
\textsf{gPC} method converts the \textsf{SSOCP} into an augmented deterministic \textsf{SOCP},  which is solved using a knotting pseudospectral
method. This method transforms the deterministic \textsf{SOCP} into an \textsf{NLP} problem, which is modeled using Pyomo \textcolor{black}{\cite{hartetal:2012:pomip}}, a Python-based software package designed to model complex optimization problems, and solved using the Interior Point OPTimizer (\textsf{IPOPT}) solver. 
In this paper, the initial guesses for the latitude, the longitude, and the heading angle are generated using the orthodromic paths between the departure and arrival locations of each flight. Typical cruise velocity and fuel consumption of the aircraft model selected are used to generate the initial guesses of the velocity and the mass of each aircraft during the flight, respectively.

\subsection{Experiment A: Three-aircraft transoceanic mission design with uncertainty in the fuel burn savings.}

Experiment A involves three transoceanic eastbound flights, Flight 1,
Flight 2, and Flight 3, with uncertainty in the fuel burn savings for the trailing aircraft. 
The three flights considered in this experiment have the following departure and arrival locations:
\begin{itemize}
\item Flight 1: New York (\texttt{JFK}) - Paris (\texttt{CDG}).
\item Flight 2: Boston (\texttt{BOS}) - Madrid (\texttt{MAD}).
\item Flight 3: Montreal (\texttt{YUL}) - London (\texttt{LHR}).
\end{itemize}
The departure times of Flight 1, Flight 2, and Flight 3 are set to 10:15, 10:30, and
10:50, respectively.
The boundary conditions for the state variables of the three aircraft
are given in Table \ref{table:boundary_conditions_flight3}. 
Flight 1, Flight 2, and Flight 3 are operated by Aircraft 1, Aircraft 2, and
Aircraft 3, respectively. 

The formation configuration is selected in advance. 
In the case of a two-aircraft formation that includes Aircraft 2,
Aircraft 2 is the leader and Aircraft 1 or
Aircraft 3 is the follower. 
In the case of a two-aircraft formation that does not include Aircraft 2, Aircraft 3 is the leader and, accordingly, Aircraft 1
is the follower.
The only three-aircraft formation allowed is the
in-line formation, where Aircraft  2 is the leader, Aircraft 3 
is the intermediate, and Aircraft 1 is the follower.

The random variable that represents the fuel burn savings for the trailing aircraft is denoted by $\theta$. This random variable is modeled as a Gaussian random variable with mean $0.1$ and standard deviation $0.02$, i.e., $\theta \sim \mathcal{N} (0.1, 0.02)$, based on several aerodynamic models and flight tests reported in \cite{kentandrichards:2021:pofffcatcs}.
The departure times are fixed.

Solving this problem entails deciding which mode of flight, i.e., formation or solo flight, is optimal,  the expected values and standard deviations of the latitude and longitude of the optimal trajectories of each aircraft as functions of time, and the expected values and standard deviation of the timing of the trajectory as functions of the orthodromic distance from the departure locations.
They include the expected values and the standard deviations of the arrival times and, in the case of formation flight, the expected values and the standard deviations of the latitude, longitude, and time of the rendezvous and splitting locations. 

The sequence of discrete states obtained in the solution is the following:


\begin{itemize}

\item
State 1: All the aircraft fly solo.

\item
State 2: Aircraft 1 and Aircraft 2 fly in a two-aircraft formation and Aircraft 3 flies solo.

\item
State 3: All the aircraft fly in a three-aircraft formation.

\item
State 4: Aircraft 1 and Aircraft 3 fly in a two-aircraft formation and Aircraft 2 flies solo.

\item
State 5: All the aircraft fly solo.

\end{itemize}

\bigskip

\begin{table}[ht!]
\centering
\caption{Experiment A: Boundary conditions of the three flights.}
\medskip
\begin{tabular}{c c c c  c}
\multicolumn{1}{ c }{\textbf{Symbol}} & \textbf{{Units}}& \textbf{{Flight 1}}             & \textbf{{Flight 2}}  & \textbf{{Flight  3}}                   \\ \hline
 $\phi_I$ & [deg]    & 40.64   & 42.36 &  45.47    \\
 $\phi_F$ & [deg]    & 48.85   & 40.48 &  51.47    \\
 $\lambda_I$ & [deg] & -73.78  & -71.06 & -73.74  \\
 $\lambda_F$ & [deg] & 2.35    & -3.57 & -0.12     \\
 $\chi_I$ & [deg]    & +54.26  &  69.25& +55.53    \\
 $V_I$ & [m/s]       & 240    &   240   & 240   \\
 $V_F$ & [m/s]       & 220   &   220    & 220   \\
 $m_I$ & [kg]        & $215 \, 000$   &  $210 \, 000$  & $220 \, 000$   \\ \hline
\end{tabular}
\label{table:boundary_conditions_flight3}
\end{table}

\begin{figure*}[ht!]
\centering
\renewcommand{\figurename}{Fig.}
{\includegraphics[width=120mm]{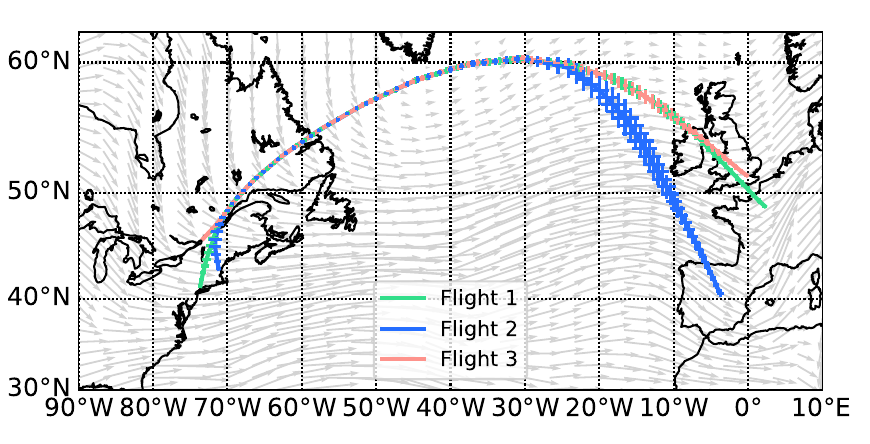}}
\caption{Experiment A: Expected values of the latitude and longitude of the optimal trajectories of each aircraft together with the corresponding 95\% confidence envelopes represented on the relevant map together with the wind field.}
\label{fig:X3_3kiwis_wind_map}
\end{figure*}

\begin{figure*}[ht!]
\centering
\renewcommand{\figurename}{Fig.}
\subfigure[Longitude]{\includegraphics[width=80mm]{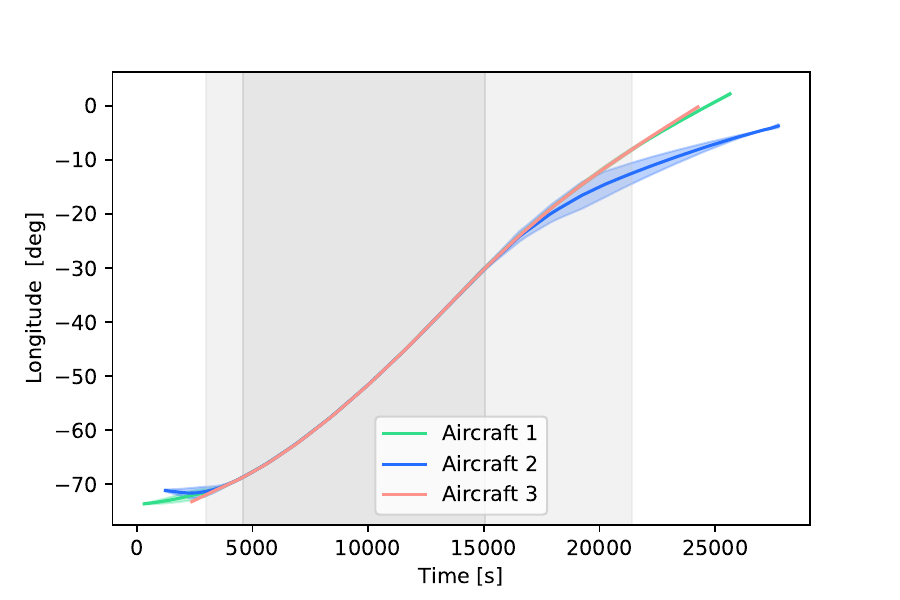}}
\subfigure[Latitude]{\includegraphics[width=80mm]{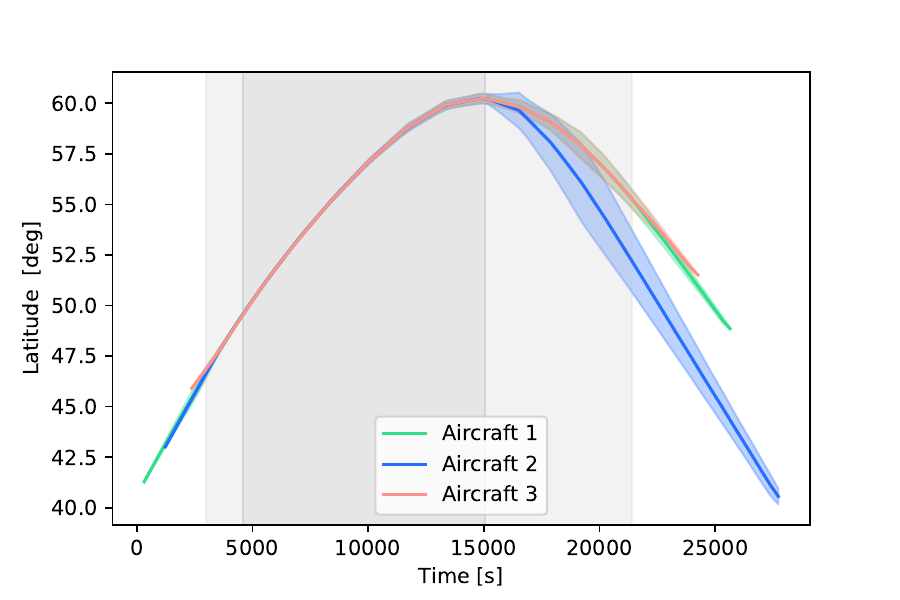}}
\subfigure[Mass]{\includegraphics[width=80mm]{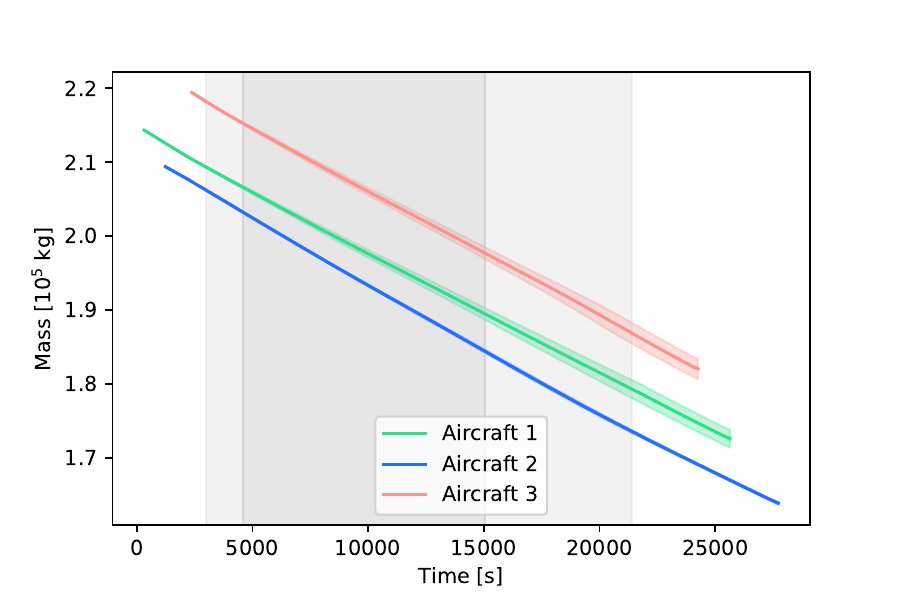}}
\subfigure[Heading angle]{\includegraphics[width=80mm]{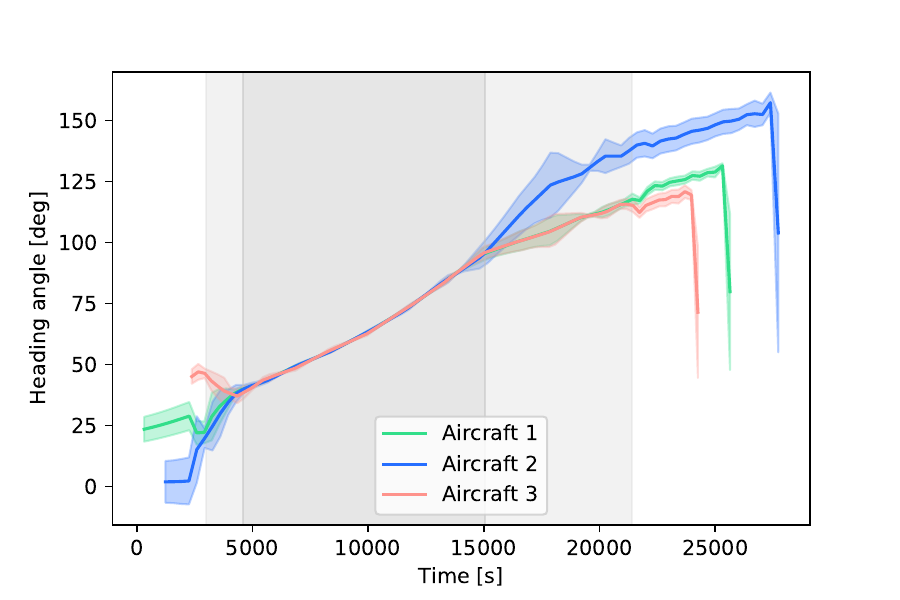}}
\subfigure[Mach number]{\includegraphics[width=80mm]{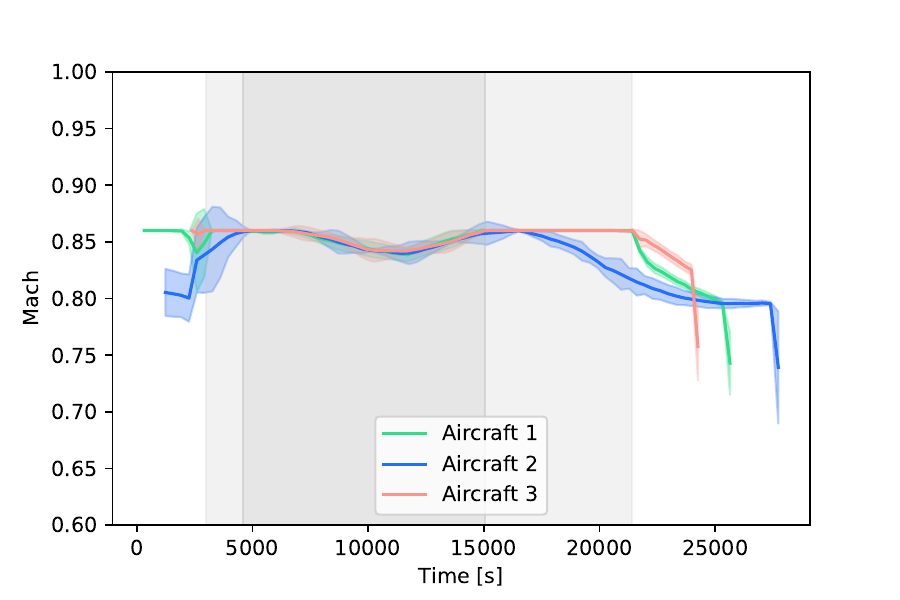}}
\caption{Experiment A: Expected values of the state variables of the optimal trajectories of each aircraft together with the corresponding 95\% confidence envelopes. }
\label{fig:state_X3_3kiwis}
\end{figure*}

\begin{figure*}[ht!]
\centering
\renewcommand{\figurename}{Fig.}
\subfigure[Lift coefficient]{\includegraphics[width=80mm]{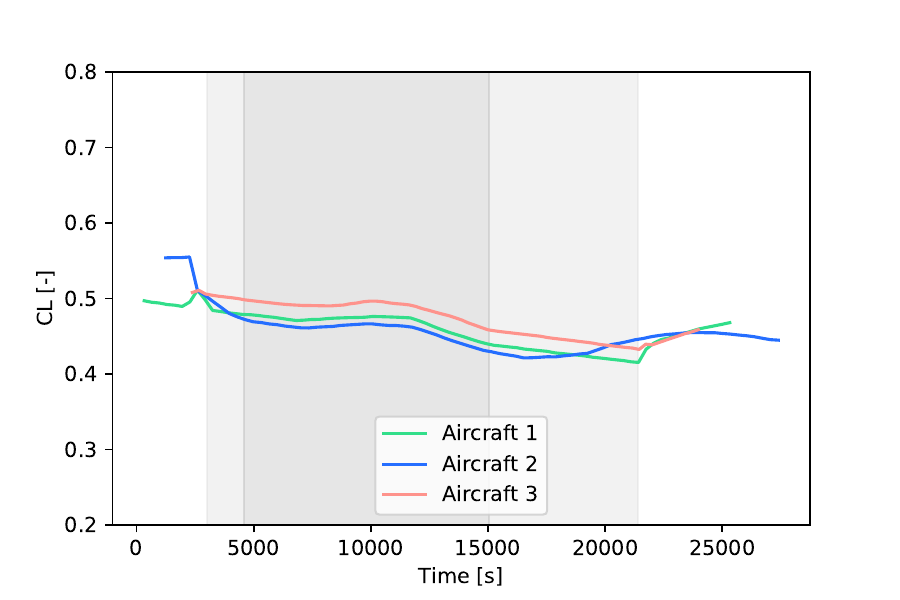}}
\subfigure[Bank angle]{\includegraphics[width=80mm]{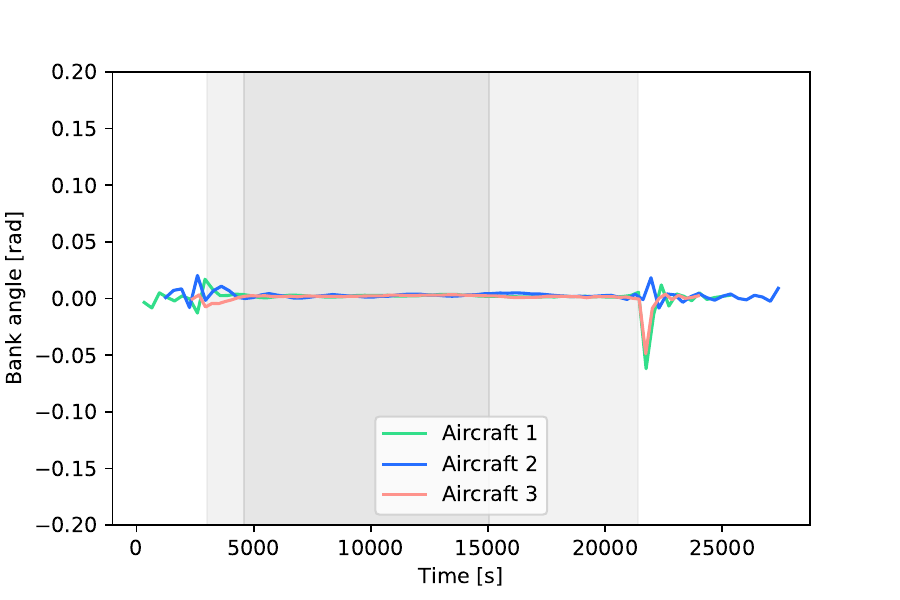}}
\subfigure[Adimensional thrust]{\includegraphics[width=80mm]{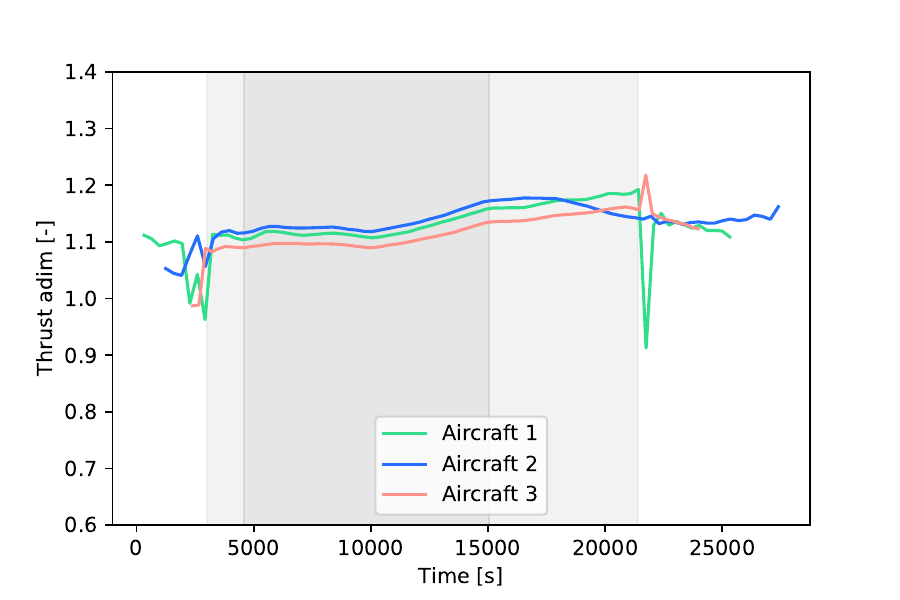}}
\caption{Experiment A:  Expected values of the control variables of each aircraft.}
\label{fig:control_X3_3kiwis}
\end{figure*}

The expected values of the  longitude and latitude of the optimal routes obtained in the solution   
together with the corresponding 95\% confidence envelopes are represented on the relevant map together with the 
wind field in Fig.~\ref{fig:X3_3kiwis_wind_map} and as functions of time in Fig.~\ref{fig:state_X3_3kiwis}.a and 
Fig.~\ref{fig:state_X3_3kiwis}.b, respectively. 
The expected values of the rest of the state variables, namely, the mass, the heading angle, and the Mach number, 
together with the corresponding 95\% confidence envelopes are represented as functions of time in
Fig.~\ref{fig:state_X3_3kiwis}.c, Fig.~\ref{fig:state_X3_3kiwis}.d, and Fig.~\ref{fig:state_X3_3kiwis}.e, respectively.
The expected values of the control variables are represented for the three flights in Fig.~\ref{fig:control_X3_3kiwis}. 
Green, blue, and red lines correspond to Flight 1, Flight 2, and Flight 3, respectively.
For a better understanding of these figures, the portions of the plots of the variables that correspond to formation flight are represented on a gray background. In particular, the portions of the plots that correspond to State 2 and State 4, which includes a two-aircraft formation, are represented on a light gray background, whereas the portion of the plots that correspond to State 3, which includes a three-aircraft formation, is represented on a dark gray background.

It is important to note that all the optimal solutions obtained considering the fuel burn savings that correspond to the nodes of the gPC expansion of the random variable $\theta \sim \mathcal{N} (0.1, 0.02)$ are formation missions with the same structure.

The expected values of the rendezvous and splitting times are reported in Table~\ref{table:rendezvous_splitting_1} together with the corresponding 95\% confidence interval.
State 2, in which Aircraft 1 and Aircraft 2 fly in formation and Aircraft 3 flies solo, has the shortest duration.
The control variables vary in a quite smooth way. 

%


Regarding the spatial variability of the solution, it can be observed in Fig.~\ref{fig:X3_3kiwis_wind_map}, Fig.~\ref{fig:state_X3_3kiwis}.a, and 
Fig.~\ref{fig:state_X3_3kiwis}.b, that 
both geographical coordinates have a similar behavior. 
In the first two discrete states, the spatial variability is rather low, and in State 3, it is almost negligible. Then, near the
first splitting location, the spatial variability starts increasing for all the flights. Finally, it decreases when the aircraft are 
approaching the final locations.
Hence, the second splitting location has a significant dependency on the random variable $\theta$ and, consequently, both the duration and the traveled distance in State 3  have a significant dependency on it. 
Flight 2  has the greatest spatial variability.

\begin{figure*}[ht!]
\centering
\renewcommand{\figurename}{Fig.}
\subfigure[Flight 1]{\includegraphics[width=80mm]{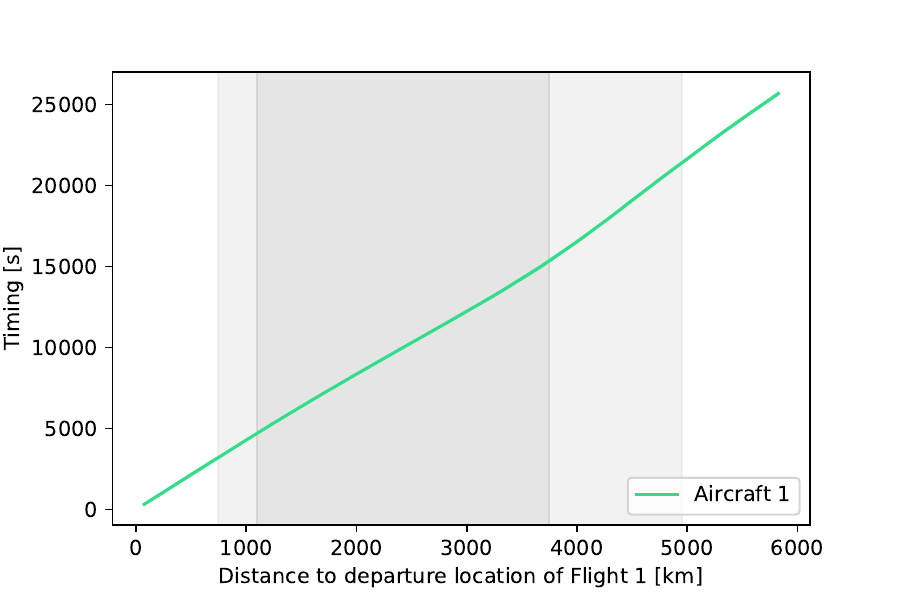}}
\subfigure[Flight 2]{\includegraphics[width=80mm]{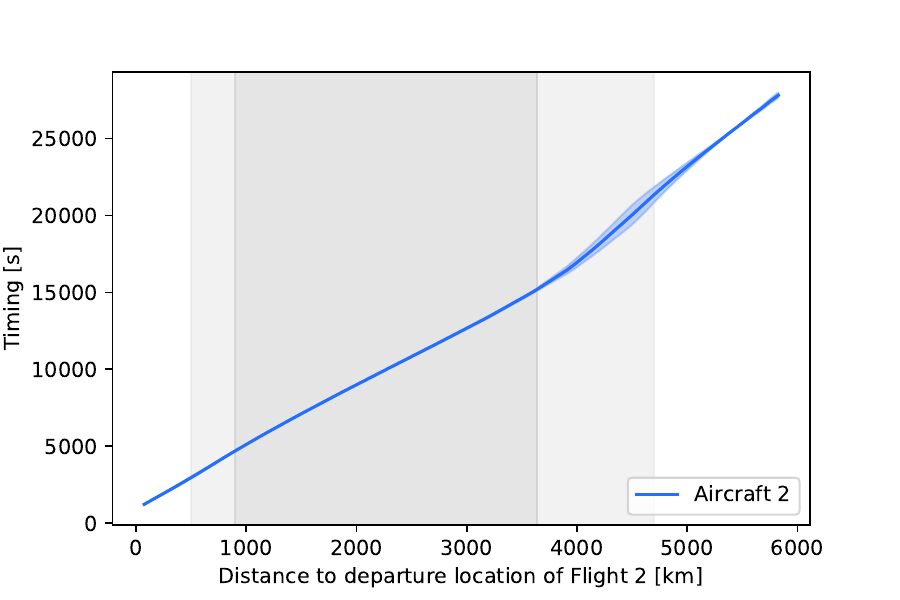}}
\subfigure[Flight 3]{\includegraphics[width=80mm]{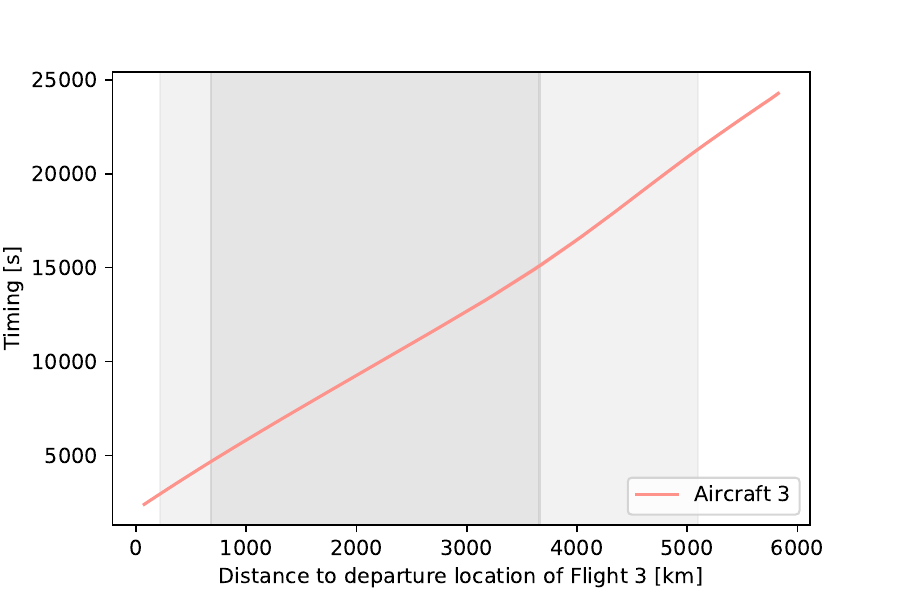}}
\caption{Experiment A: Expected values of the timing of the optimal trajectories of each aircraft together with the corresponding 95\% envelopes.}
\label{fig:time_deviation_X3_3kiwis}
\end{figure*}

In addition to the spatial variability, the temporal variability is also quantified in order to have complete information regarding the spatio-temporal variability of the solution. The expected values of the timing of the solution together with the corresponding 95\% confidence envelopes are represented in Fig.~\ref{fig:time_deviation_X3_3kiwis} as functions of the orthodromic distance from the departure location of each flight. It can be observed in this figure that the temporal variability is very low in all the flights. In particular, in Flight 1 and Flight 3, the temporal variability is almost negligible during the whole flight time, whereas, in Flight 2, it starts increasing after the first splitting location. However, despite this increase, it is very low in Flight 2 as well, leading to the
conclusion that the temporal variability of all the flights slightly depends on the random variable $\theta$.

It can be observed in Table~ \ref{table:rendezvous_splitting_1}, in which the expected values and the 95\% confidence intervals of the rendezvous and splitting times are reported, that the amplitudes of these intervals are lower for the rendezvous times than for the splitting times, but remains fairly small for both of them.

\begin{table}[h!]
\caption{Experiment A: Expected values and 95\% confidence intervals of the rendezvous and splitting times.}
\centering
\begin{tabular}{lccc}
  & \textbf{Expected} & \textbf{ 95\% confidence} &   \\
    & \textbf{value} & \textbf{ interval} &   \\
  \hline
\textbf{First rendezvous time [h]} &  0.88 &  [0.87, 0.89]   \\ 
\textbf{Second rendezvous time [h]} &  1.32 &  [1.32, 1.32]   \\ 
 \textbf{First splitting time [h]} &  4.19 &   [4.17, 4.21]   \\ 
\textbf{Second splitting time [h]} &  5.95 &  [5.92, 5.98]  
\end{tabular}
\label{table:rendezvous_splitting_1}
\end{table}

The expected values and the 95\% confidence intervals for both the flight time, expressed in hours, and the fuel burn, expressed in tonnes, for each flight are listed in Table~\ref{table:results_comparison_3kiwis_X3}. 
From the amplitudes of the 95\% confidence intervals, it is easy to see that, as already mentioned, the flight time variability is rather low. In contrast, as expected, the uncertainty in the fuel burn savings has a significant impact on the  variability of the fuel consumption.

\begin{table*}[ht!]
\centering
\caption{Experiment A: Expected values and 95\% confidence intervals of flight times and fuel consumptions of each flight.}
\medskip
\begin{tabular}{c ccc cc}
\multicolumn{1}{c }{}            & \multicolumn{2}{c }{\textbf{Flight time} \textbf{{[}}$\boldsymbol{h}$\textbf{{]}}}  & \multicolumn{2}{c }{\textbf{Fuel burn} \textbf{{[}}$\boldsymbol{t}$\textbf{{]}} }    \\  
\cline{2-5}
\multicolumn{1}{ c }{\textbf{ }} & \textbf{Expected value} & \textbf{95\% confidence interval}  & \textbf{Expected value} & \textbf{95\% confidence interval}  &    \\ \hline
\textbf{Flight 1}    &   7.13  &  [7.12, 7.15]    &     42.47  &   [41.19, 43.74]  \\
\textbf{Flight 2}    &    7.47  &  [7.43, 7.52]  &     46.24  &   [45.90, 46.58]   \\
\textbf{Flight 3}    &     6.16  &  [6.15, 6.18]   &    38.00  &   [36.58, 39.42]  \\
  \hline
\end{tabular}
\label{table:results_comparison_3kiwis_X3}
\end{table*}


To quantify the benefits of formation flight with respect to solo flight, the expected values of the flight time and the fuel consumption, along with the expected value of the corresponding \textsf{DOC} expressed in monetary units, \textit{mu}, are estimated assuming that the three flights are performed as
solo flights. The obtained results are reported in Table~\ref{table:DOC_solo_FF}. For the sake of comparison, the expected values of the \textsf{DOC} 
obtained in both formation and solo flights of the three aircraft are reported in Table~\ref{table:results_comparison_FF_SF_3kiwis}. It is remarkable that, in spite of
the uncertainties considered in the fuel burn savings for the trailing aircraft, the reduction in the total \textsf{DOC}  achieved with
formation flight amounts to $2.16\%$.

\begin{table}[h!]
\caption{Experiment A: Expected values of flight times, fuel consumptions, and \textsf{DOC} in solo flight of the three aircraft.}
\centering
\begin{tabular}{lccc}
  & {\textbf{Flight time} \textbf{{[}}$\boldsymbol{h}$\textbf{{]}}} & 
 {\textbf{Fuel burn} \textbf{{[}}$\boldsymbol{t}$\textbf{{]}}} &   
 {\textbf{DOC} \textbf{{[}}$\boldsymbol{mu}$\textbf{{]}}}\\
  \hline
\textbf{Flight 1} &  7.06 &  45.73 &  39635.79  \\ 
\textbf{Flight 2} &  7.32 &  45.44 &  39713.60   \\ 
\textbf{Flight 3} &  6.02 &  39.52 &  34165.60
\end{tabular}
\label{table:DOC_solo_FF}
\end{table}

\begin{table*}[ht!]
\centering
\caption{Experiment A: Expected values of the  \textsf{DOC} in solo and formation flights of the three aircraft.}
\medskip
\begin{tabular}{l llllll}
\multicolumn{1}{ c }{\textbf{ }} & \textbf{Flight 1} & \textbf{Flight 2}  & \textbf{Flight 3} & \textbf{Total}  \\ 
\cline{1-5}
{\textbf{DOC  solo flight} \textbf{{[}}$\boldsymbol{mu}$\textbf{{]}}}  &   39635.79   & 39713.60    &     34165.60  &    113514.99
 \\
{\textbf{DOC formation flight} \textbf{{[}}$\boldsymbol{mu}$\textbf{{]}}}  &   37429.40  & 40435.59   &     33252.80  &   111117.79
 \\
{$\boldsymbol{\Delta}$\textbf{DOC} \textbf{{[}}$\boldsymbol{\%}$\textbf{{]}}}   &   -5.89 &  +1.79    &     -2.75    &   -2.16
\\
\cline{1-5}
\end{tabular}
\label{table:results_comparison_FF_SF_3kiwis}
\end{table*}

To quantify the effects of the uncertainty in the fuel burn savings for the trailing aircraft on the \textsf{DOC}, the same formation mission  is designed in the absence of uncertainty, in which the fuel burn savings for both the trailing and intermediate aircraft are set to 10\%. 
For the sake of comparison, the obtained results in terms of the \textsf{DOC} for the deterministic and stochastic formation missions are summarized in Table~\ref{table:results_comparison_FF_deterministic_3kiwis}.  
It can be seen that the increment of the total \textsf{DOC} due to the presence of uncertainties amounts to 1.60\%.

\begin{table*}[ht!]
\centering
\caption{Experiment A:  Values of the  \textsf{DOC} in deterministic and stochastic formation missions of the three aircraft.}
\medskip
\begin{tabular}{l llllll}
\multicolumn{1}{ c }{\textbf{ }} & \textbf{Flight 1} & \textbf{Flight 2}  & \textbf{Flight 3} & \textbf{Total}  \\ 
\cline{1-5}
{\textbf{DOC  deterministic formation mission} \textbf{{[}}$\boldsymbol{mu}$\textbf{{]}}}  &   36178.9 & 40205.67 & 32980.96  &    109365.53
 \\
{\textbf{DOC stochastic formation mission} \textbf{{[}}$\boldsymbol{mu}$\textbf{{]} (expected values)} } &   37429.40  & 40435.59   &     33252.80  &   111117.79
 \\
{$\boldsymbol{\Delta}$\textbf{DOC} \textbf{{[}}$\boldsymbol{\%}$\textbf{{]}}}   &   +3.45 &  +0.57    &     +0.82    &   +1.60
\\
\cline{1-5}
\end{tabular}
\label{table:results_comparison_FF_deterministic_3kiwis}
\end{table*}


\subsection{Experiment B: Two-aircraft transoceanic mission design with uncertainty in \textcolor{black}{the} departure times of the flights.}

Experiment B involves two transoceanic eastbound flights, Flight 1 and Flight 2, with uncertainty in the departure times of the flights. The flights considered in this experiment have the
same departure and arrival locations and the same boundary values for the
state variables as Flight 1 and Flight 2 in Experiment A. They have the following departure and arrival locations:
\begin{itemize}
\item Flight 1: New York (\texttt{JFK}) - Paris (\texttt{CDG}).
\item Flight 2: Boston (\texttt{BOS}) - Madrid (\texttt{MAD}).
\end{itemize}
The departure times of Flight 1 and Flight 2 are set to 10:15 and 10:30, respectively.
The boundary conditions for the state variables of the two aircraft
are given in Table \ref{table:boundary_conditions_flight3}. 
Flight 1 and Flight 2 are operated by Aircraft 1 and Aircraft 2, respectively.

As in Experiment A, the formation configuration is selected in advance. 
In the case of formation flight, Aircraft 1 is the trailing
aircraft and Aircraft 2 is the leader. 
In this experiment, the fuel burn savings for the trailing aircraft are set to 10\%.

The random variables that represent the delays with respect to the scheduled departure times of Flight 1 and Flight 2 from \texttt{JFK} and \texttt{BOS} airports are denoted by $\theta_1$ and $\theta_2$, respectively. These random variables have been modeled using a Gaussian mixture distribution, the parameters of which are estimated from real delay data obtained
from the {\sf Transtats} online database\footnote{\href{https://www.transtats.bts.gov/
}{https://www.transtats.bts.gov/
}} of the U.S. Bureau of Transportation Statistics. This database provides detailed information on the
U.S. transportation systems, including data on departure delays by airport and airline.
The expectation-maximization algorithm for fitting mixture-of-Gaussian models to data is employed for this purpose \cite{murphy:2012:mlapp}. Departure delay data corresponding to January 2020 are used. 





\begin{figure*}[ht!]
\centering
\renewcommand{\figurename}{Fig.}
\subfigure[$\theta_2$, departure delay of Flight 1 from \texttt{JFK} airport]{\includegraphics[width=81mm]{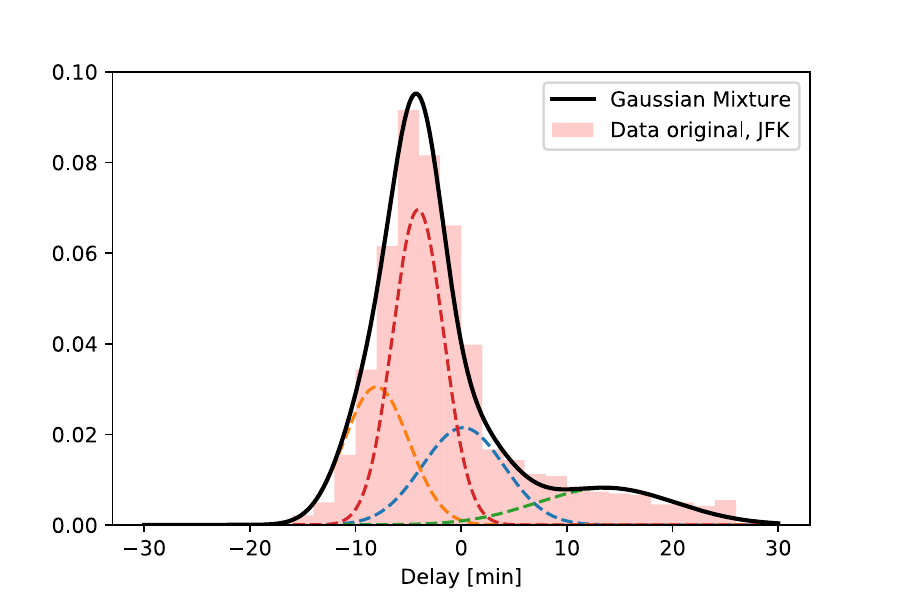}}
\subfigure[$\theta_2$, departure delay of Flight 2 from \texttt{BOS}  airport]{\includegraphics[width=81mm]{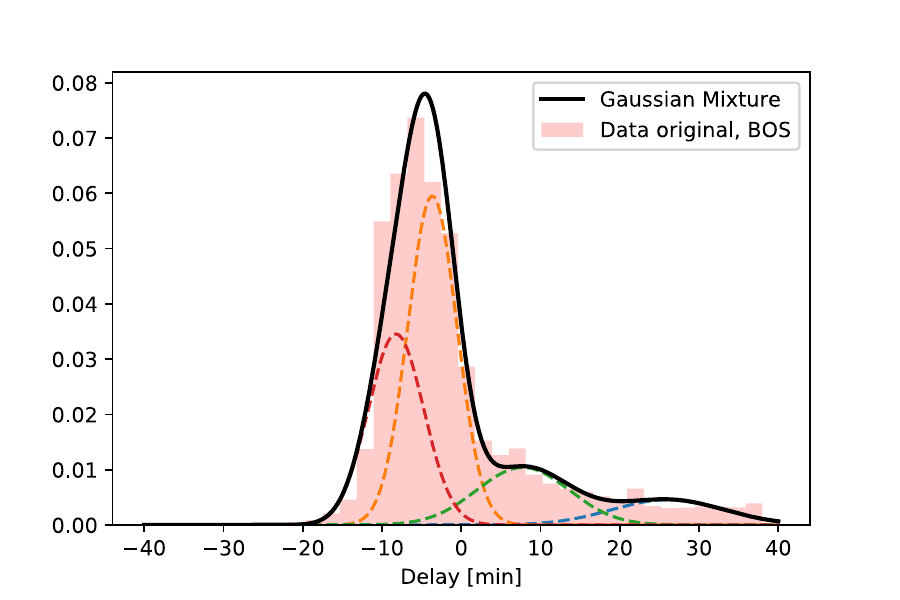}}
\caption{Experiment B: Estimated probability density functions of the random variables $\theta_1$ and $\theta_2$ that represent the departure delays of each flight. }
\label{fig:estimated_distrib}
\end{figure*}

The resulting Gaussian mixture probability density functions that model the random variables $\theta_1$ 
and $\theta_2$ that represent the departure delays at \texttt{JFK} and \texttt{BOS} airports both have 4 components. The values of the parameters of the Gaussian mixture probability density functions are given in Table~\ref{table: GMM_parameters}, where $p_{i,j}, i=1,2, j=1,2,3,4$ are the weights of the Gaussian component probability density functions of $\theta_i$, which are probabilities that sum to 1, and $\mu_{i,j}$ and $\sigma_{i,j}$ are their means and standard deviations, respectively. This result is consistent with previous studies on estimation of flight departure delay distributions \cite{tu2008estimating}.
Fig.~\ref{fig:estimated_distrib} shows in black the obtained Gaussian mixture distribution of the departure delays of  Flight 1 and Flight 2 from \texttt{JFK} and \texttt{BOS} airports, respectively, 
together with the histogram of the data. A bin width of 2 minutes is used for the histograms. 
The components of the Gaussian mixture distribution are represented by dotted lines. 
The mean and the standard deviation of the random variable $\theta_1$ are  -1.67 [min]    and  7.69 [min], respectively.
The mean and the standard deviation of the random variable  $\theta_2$ are  -0.88  [min] and 10.37 [min], respectively.

\begin{table*}[ht!]
\centering
\caption{Experiment  B: Estimated parameters of the Gaussian mixture probability density functions
 that model the departure delays of the flights.  }
\medskip
\begin{tabular}{c c c ccc cc}
& \multicolumn{1}{ c } {\textbf{Departure airport}} & \textbf{Random variable} & $\boldsymbol{p_{i,1}, p_{i,2}, p_{i,3}, p_{i,4}}$  & $\boldsymbol{\mu_{i,1}, \mu_{i,2}, \mu_{i,3}, \mu_{i,4}}$ & $\boldsymbol{\sigma_{i,1}, \sigma_{i,2}, \sigma_{i,3}, \sigma_{i,4}}$  &    \\ \hline

\textbf{Flight 1}& \texttt{JFK}    & $\theta_1$   & 0.39, 0.17, 0.27, 0.17 &  -4.94, 11.94, -0.99, -8.91 &     2.20, 7.17, 2.93, 2.89 \\
\textbf{Flight 2}&{\texttt{BOS}}    &$\theta_2$  & 0.24, 0.56, 0.06,  0.14 & -1.76,  -6.61, 28.67, 10.92  &     3.39, 3.70, 5.68, 5.78  \\
  \hline
\end{tabular}
\label{table: GMM_parameters}
\end{table*}


\bigskip

\bigskip

As in Experiment A, solving this problem entails deciding which mode of flight is optimal, the expected values and standard deviations of the latitude and longitude of the optimal trajectories of each aircraft as functions of time, and the expected values and standard deviation of the timing of the trajectory as functions of the distance.

The sequence of discrete states obtained in the solution is the following:

\begin{itemize}
\item
State 1: Both aircraft fly solo.
\item
State 2: The aircraft fly in formation.
\item
State 3: Both aircraft fly solo.
\end{itemize}



\begin{figure*}[ht!]
\centering
\renewcommand{\figurename}{Fig.}
{\includegraphics[width=120mm]{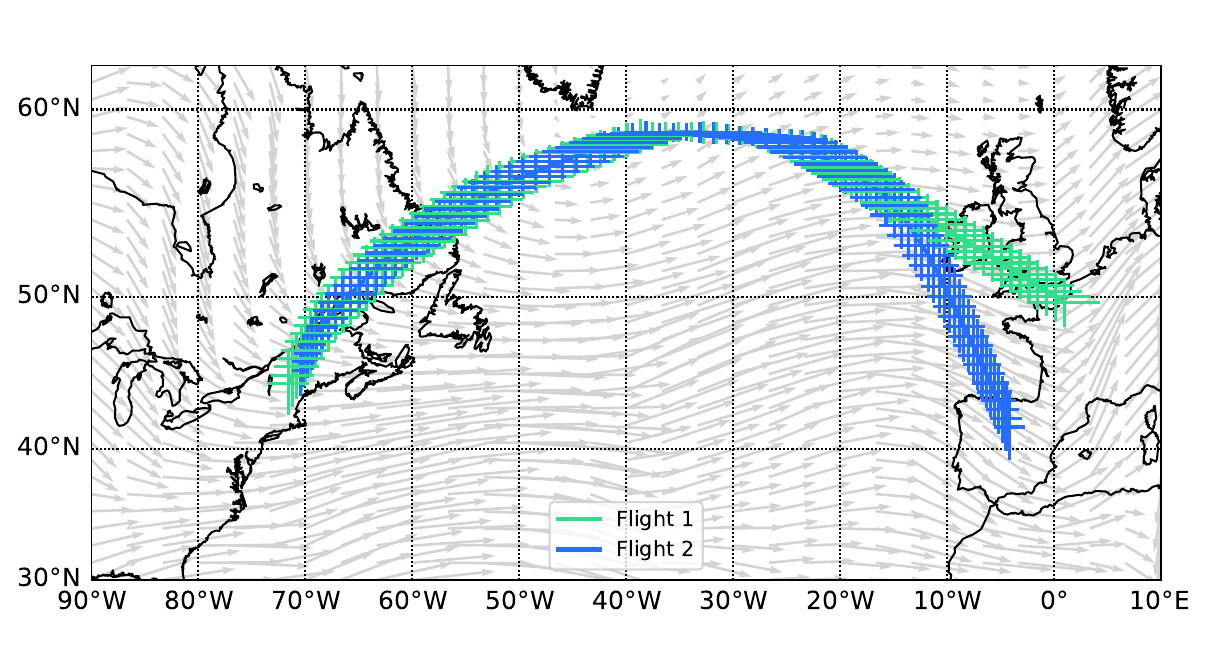}}
\caption{Experiment B: Expected values of the longitude and latitude of the optimal trajectories  of each aircraft together with the corresponding 95\% confidence envelopes represented on the relevant map together with the wind field.} 
\label{fig:X9_2kiwis_wind_map}
\end{figure*}

\begin{figure*}[ht!]
\centering
\renewcommand{\figurename}{Fig.}
\subfigure[Longitude]{\includegraphics[width=80mm]{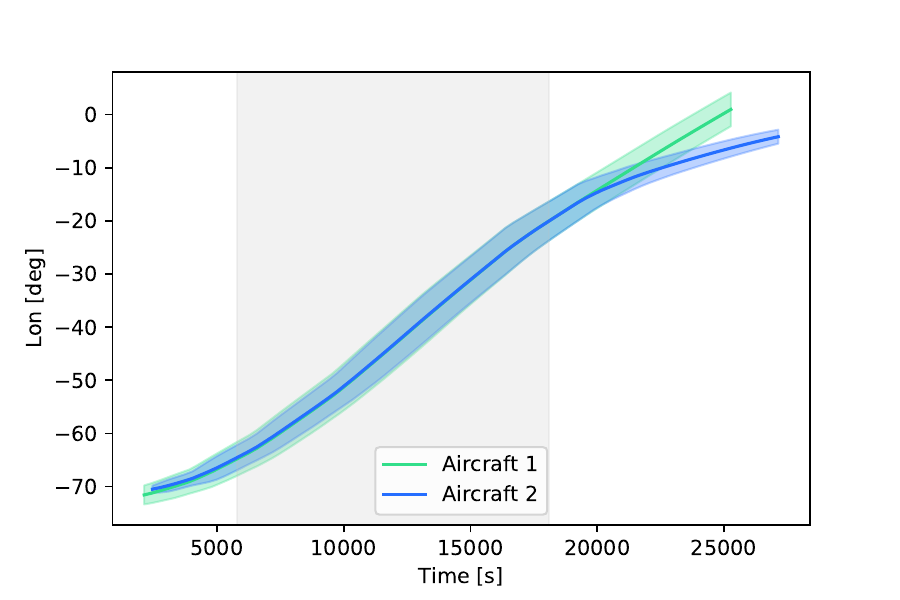}}
\subfigure[Latitude]{\includegraphics[width=80mm]{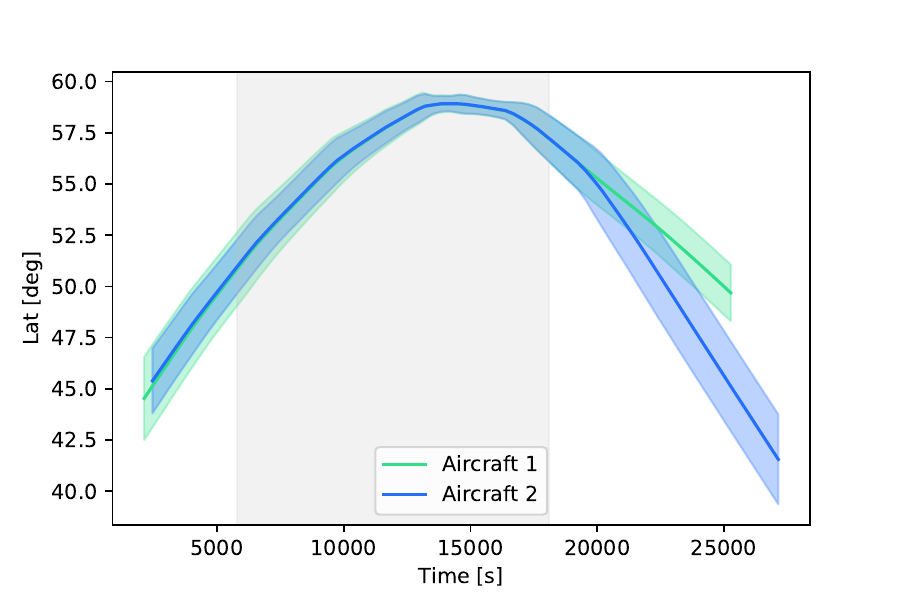}}
\subfigure[Mass]{\includegraphics[width=80mm]{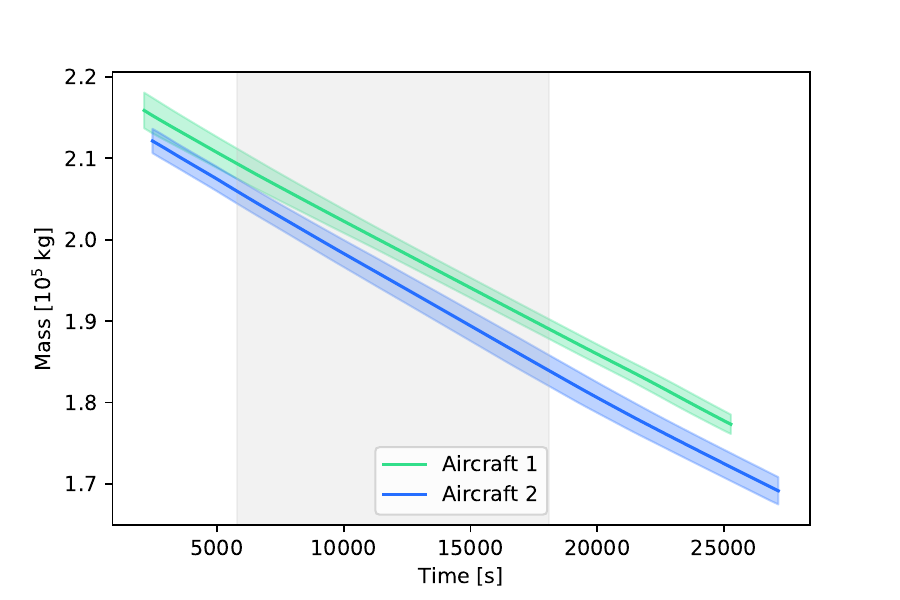}}
\subfigure[Heading angle]{\includegraphics[width=80mm]{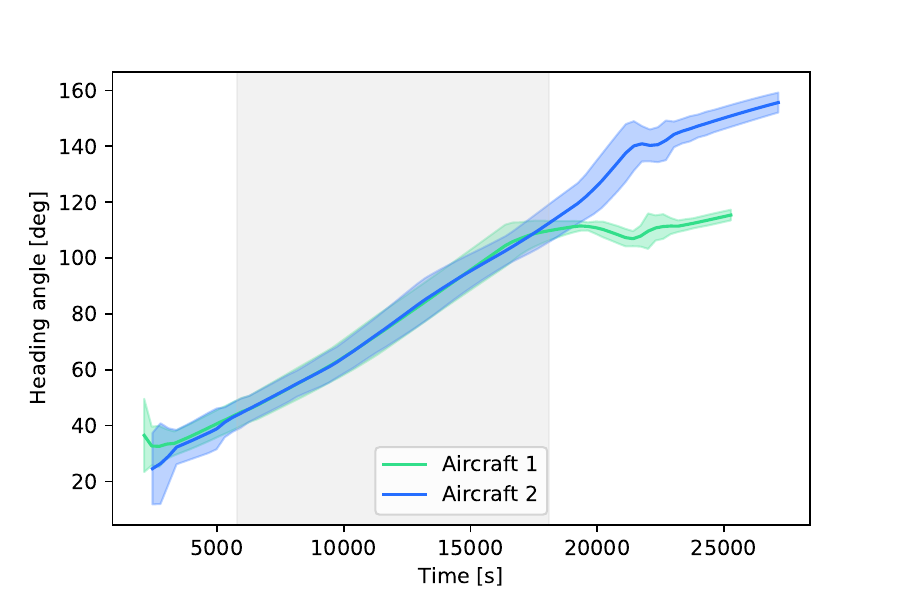}}
\subfigure[Mach number]{\includegraphics[width=80mm]{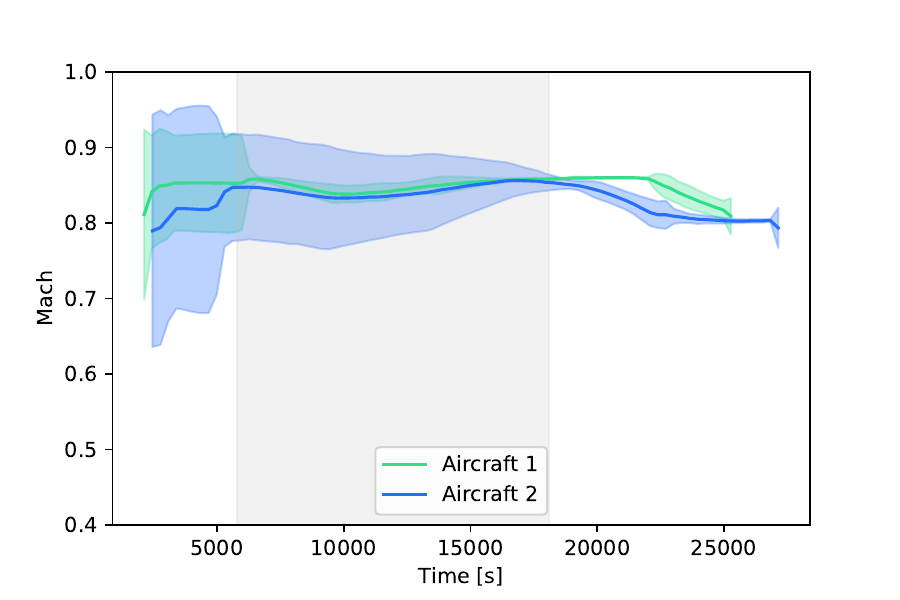}}
\caption{Experiment B: Expected values of the state variables of the optimal trajectories of each aircraft together with the corresponding 95\% confidence envelopes.}
\label{fig:state_X9_2kiwis}
\end{figure*}

\begin{figure*}[ht!]
\centering
\renewcommand{\figurename}{Fig.}
\subfigure[Lift coefficient]{\includegraphics[width=80mm]{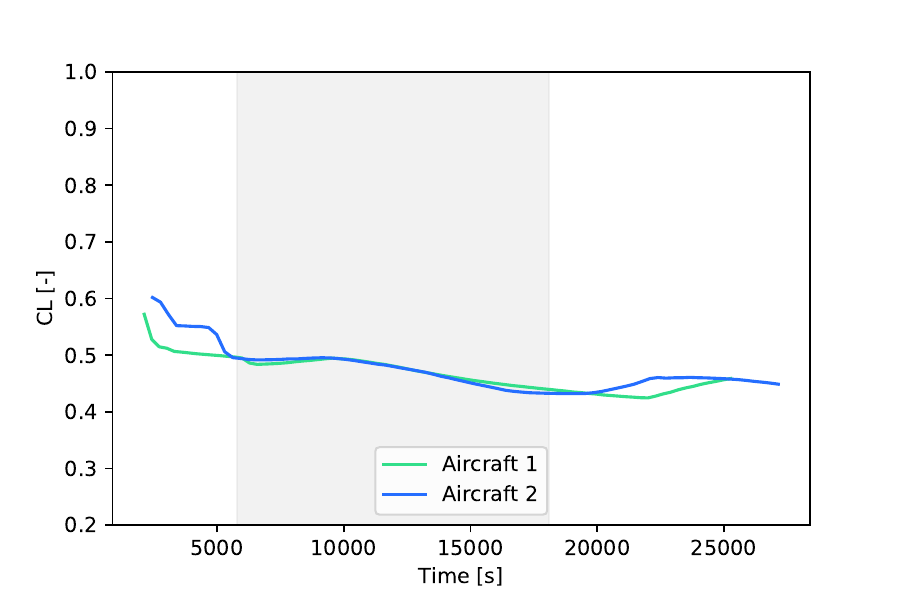}}
\subfigure[Bank angle]{\includegraphics[width=80mm]{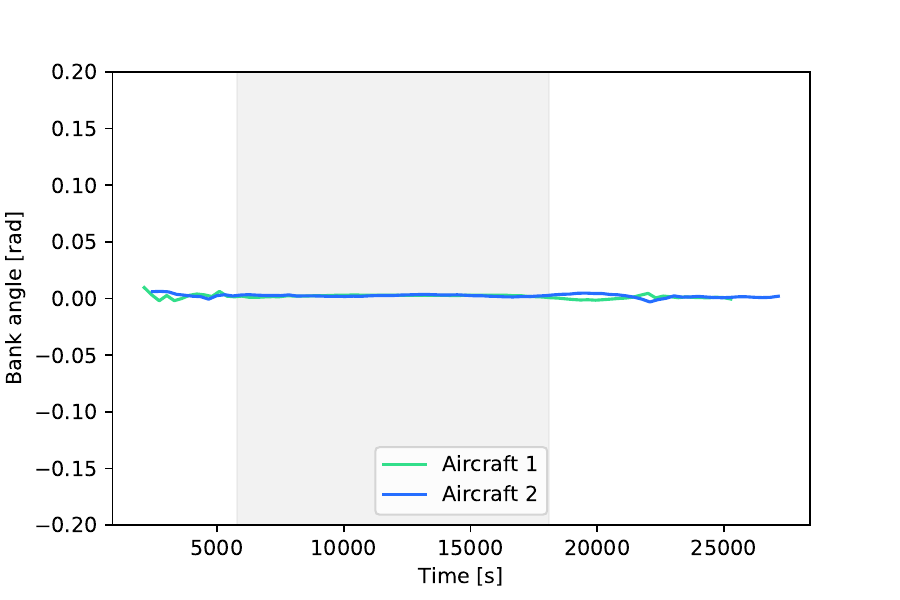}}
\subfigure[Adimensional thrust]{\includegraphics[width=80mm]{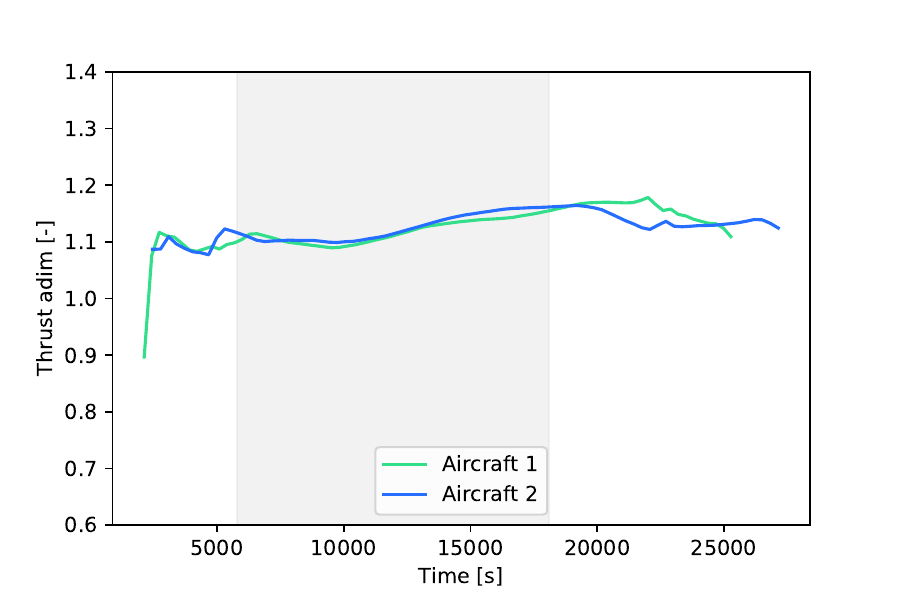}}
\caption{Experiment B:  Expected values of the control variables of each aircraft.}
\label{fig:control_X9_2kiwis}
\end{figure*}

The expected values of the  longitude and latitude of the optimal routes obtained in the solution together with the corresponding 95\% confidence envelopes are represented on the relevant map together with the wind field in 
Fig.~\ref{fig:X9_2kiwis_wind_map} and as functions of time in Fig.~\ref{fig:state_X9_2kiwis}.a and Fig.~\ref{fig:state_X9_2kiwis}.b, respectively.
The expected values of the rest of the state variables, namely, the mass, the heading angle, and the Mach number, together with the corresponding 95\% confidence envelopes are represented as functions of time in 
Fig.~\ref{fig:state_X9_2kiwis}.c, Fig.~\ref{fig:state_X9_2kiwis}.d, and Fig.~\ref{fig:state_X9_2kiwis}.e, respectively.
The expected values of the control variables are represented for the two flights in Fig.~\ref{fig:control_X9_2kiwis}. 
Green and blue lines correspond to Flight 1 and Flight 2, respectively. 
As in Experiment A, the portions of the plots of the variables that correspond to State 2, in which the aircraft fly in formation, are represented on a gray background.

It is important to note that all the optimal solutions obtained considering the fuel burn savings that correspond to the nodes of the gPC expansion of the two random variables are formation missions with the same structure.

The expected values of the rendezvous and splitting times are reported in Table~\ref{table:rendezvous_splitting_2} together with the corresponding 95\% confidence intervals. State 2, in which the aircraft fly in formation, has the longest duration.


The control variables vary in a quite smooth way.
Comparing the expected values of the state variables and the corresponding 95\% confidence envelopes represented as functions of time in Fig.~\ref{fig:state_X9_2kiwis}
with the corresponding ones of Experiment A given in Fig.~\ref{fig:state_X3_3kiwis}, 
the high variability in
every state variable of each flight is striking. In particular, the
Mach number variability in the early part of the flight is especially
worthy of note. This is due to the fact that anticipations and delays in the departure times of the flights
are compensated by increments and reductions of their velocities.


Regarding the spatial variability of the solution, it can be observed in 
Fig.~\ref{fig:X9_2kiwis_wind_map}, 
Fig.~\ref{fig:state_X9_2kiwis}.a, and 
Fig.~\ref{fig:state_X9_2kiwis}.b
that the longitude variability is greater in State 2, in which the aircraft fly in formation. 
In contrast, the latitude variability has a minimum in State 2, which is located around the maximum value of this geographical coordinate. Aircraft 2, the leading aircraft, has the greatest spatial variability.

\begin{figure*}[ht!]
\centering
\renewcommand{\figurename}{Fig.}
\subfigure[Flight 1]{\includegraphics[width=80mm]{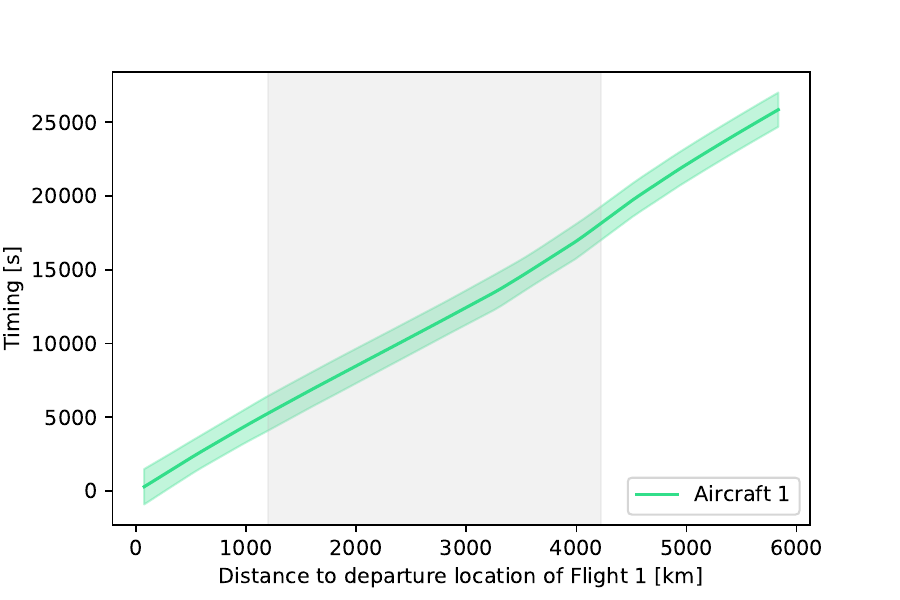}}
\subfigure[Flight 2]{\includegraphics[width=80mm]{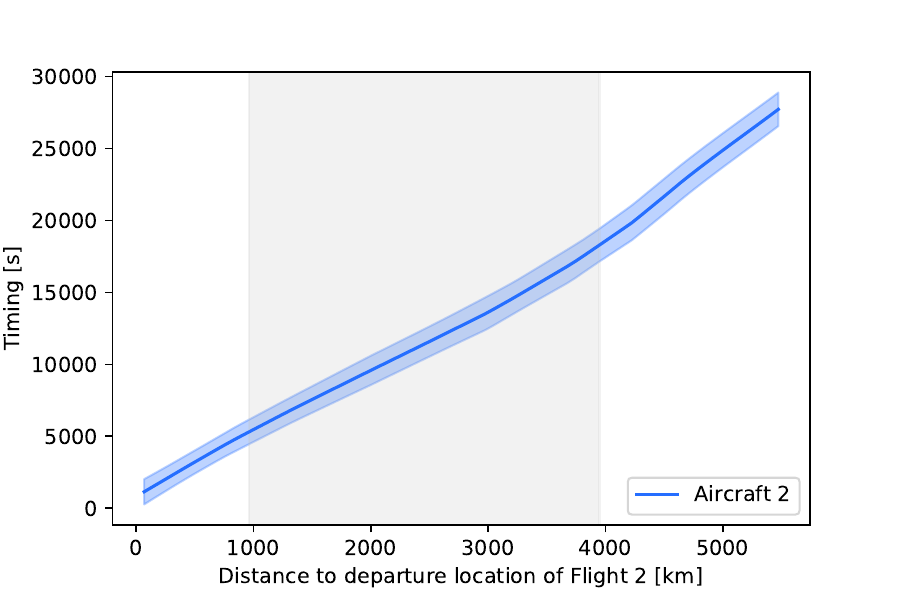}}
\caption{Experiment B: Expected values of the timing of the optimal trajectories of each aircraft together with the corresponding 95\% envelopes.} 
\label{fig:time_deviation_X9_2kiwis}
\end{figure*}

As in Experiment A, 
in addition to the spatial variability, the temporal variability is also quantified in order to have complete information regarding the spatio-temporal variability of the solution.
The expected values of the timing of the solution together with the corresponding 95\% confidence envelopes are represented in Fig.~\ref{fig:time_deviation_X9_2kiwis}, as functions of the orthodromic distance from the departure locations of each flight. 
It can be observed in this figure that the temporal variability has a similar behavior for both flights, remaining nearly
constant throughout the flight. As in the case of spatial variability, the temporal variability of the flights in Experiment B is considerably
greater than the temporal variability of the flights in Experiment A. As expected, uncertainty in the departure times of the
flights has a significant impact on the temporal variability of the trajectories, more than the uncertainty in the fuel burn savings
for the trailing aircraft.

It can be observed in Table~\ref{table:rendezvous_splitting_2}, in which the expected values and the 95\% confidence intervals of the rendezvous and splitting times are reported, that the amplitudes of these intervals are similar for the rendezvous and splitting times, being, in both cases, considerably greater than the amplitudes of the \textcolor{black}{95\% confidence} intervals of the rendezvous and splitting times obtained in Experiment A.

\begin{table}[h!]
\caption{{Experiment B: Expected values and 95\% confidence intervals of the rendezvous and splitting times. }}
\centering
\begin{tabular}{llll}
  & \textbf{Expected value} & \textbf{95\% confidence interval}   \\
  \hline
\textbf{Rendezvous time [h]} &  1.61 &  [1.29, 1.94]   \\
\textbf{Splitting time [h]} & 5.03 &  [4.71, 5.34]   \\
\end{tabular}
\label{table:rendezvous_splitting_2}
\end{table}

The expected values and the 95\% confidence intervals for both the flight time, expressed in hours, and the fuel burn, expressed in tonnes, for each flight are listed in Table~\ref{table: results_comparison_X9_2kiwis}. From the amplitudes of the 95\% confidence intervals, it is easy to see that, as already mentioned, the flight time variability is rather high. 
Additionally, the uncertainty in the departure times also has a significant impact on the variability of the fuel consumption.

\begin{table*}[ht!]
\centering
\caption{Experiment B: Expected values and 95\% confidence intervals of flight times and fuel consumptions of each flight.}
\medskip
\begin{tabular}{c ccc cc}
\multicolumn{1}{c }{}            & \multicolumn{2}{c }{\textbf{Flight time} \textbf{{[}}$\boldsymbol{h}$\textbf{{]}}}  & \multicolumn{2}{c }{\textbf{Fuel burn} \textbf{{[}}$\boldsymbol{t}$\textbf{{]}} }    \\  
\cline{2-5}
\multicolumn{1}{ c }{\textbf{ }} & \textbf{Expected value} & \textbf{95\% confidence interval}  & \textbf{Expected value} & \textbf{95\% confidence interval}  &    \\ \hline
\textbf{Flight 1}    &   7.18  &  [6.86, 7.50]    &    43.58  &   [42.52, 44.63]  \\
\textbf{Flight 2}    &    7.70 &  [7.37, 8.02]  &     46.69 &   [46.27, 47.11]   \\
  \hline
\end{tabular}
\label{table: results_comparison_X9_2kiwis}
\end{table*}

To quantify the benefits of formation flight with respect to solo flight, the expected values of the flight time and the fuel consumption along with the expected value of the
corresponding \textsf{DOC} expressed in monetary units, \textit{mu}, for the two solo flights, given in the first two rows of Table~\ref{table:DOC_solo_FF}, are compared to the expected values of the \textsf{DOC}
obtained in formation flight. The obtained results are reported in Table ~\ref{table:results_comparison_FF_SF_2kiwis}. It can be observed that, in the presence
of uncertainties in the departure times of the flights, 
the probability density funcion of which is estimated from real departure delay data, small benefits are expected in terms of the reduction of the total \textsf{DOC},
which amount to only 0.11\%. However, in real-life scenarios, some mitigation measures could be implemented to reduce the
negative impact of flight delays on the benefits of formation flight, such as adjusting the flight plan of the formation mission
before departure when one of the flights of a formation is delayed.

\begin{table}[ht!]
\centering
\caption{Experiment B: Expected values of the \textsf{DOC} in solo and formation flights of the two aircraft.}
\medskip
\begin{tabular}{l llllll}
\multicolumn{1}{ c }{\textbf{ }} & \textbf{Flight 1} & \textbf{Flight 2}  & \textbf{Total}  \\ 
\cline{1-4}
{\textbf{DOC  solo flight} \textbf{{[}}$\boldsymbol{mu}$\textbf{{]}}}  &   39635.79  & 39713.60    &     79349.39
 \\
{\textbf{DOC formation flight} \textbf{{[}}$\boldsymbol{mu}$\textbf{{]}}}  &   38260.39  & 40999.01  &   79259.39
 \\
{$\boldsymbol{\Delta}$\textbf{DOC} \textbf{{[}}$\boldsymbol{\%}$\textbf{{]}}}   &  -3.59   &     +3.14\%  &   -0.11 
\\
\cline{1-4}
\end{tabular}
\label{table:results_comparison_FF_SF_2kiwis}
\end{table}

To quantify the effects of the uncertainty in the departure delays on the \textsf{DOC}, the same formation mission  is designed in the absence of uncertainty, in which no delays are considered. 
For the sake of comparison, the obtained results in terms of the \textsf{DOC} for the deterministic and stochastic formation missions are summarized in Table~\ref{table:results_comparison_FF_deterministic_2kiwis}.  
It can be seen that the increment of the total \textsf{DOC} due to the presence of uncertainties amounts to 3.76\%.


\begin{table*}[ht!]
\centering
\caption{Experiment B:  Values of the  \textsf{DOC} in deterministic and stochastic formation missions of the two aircraft.}
\medskip
\begin{tabular}{l llllll}
\multicolumn{1}{ c }{\textbf{ }} & \textbf{Flight 1} & \textbf{Flight 2}  & \textbf{Total}  \\ 
\cline{1-4}
{\textbf{DOC  deterministic formation flight} \textbf{{[}}$\boldsymbol{mu}$\textbf{{]}}}  &   36178.90 & 40205.67 &    76384.57
 \\
{\textbf{DOC  stochastic formation flight} \textbf{{[}}$\boldsymbol{mu}$\textbf{{]} (expected values)}}  &   38260.39  & 40999.01  &   79259.39
 \\
{$\boldsymbol{\Delta}$\textbf{DOC} \textbf{{[}}$\boldsymbol{\%}$\textbf{{]}}}   &   +5.75 &  +1.97        &   +3.76
\\
\cline{1-4}
\end{tabular}
\label{table:results_comparison_FF_deterministic_2kiwis}
\end{table*}

Based on the results of the two experiments, it is possible to conclude that uncertainties have a significant impact on the formation flight benefits in terms of the reduction of the \textsf{DOC}. In particular, the uncertainties in the departure times have a far greater impact on the \textsf{DOC} than the uncertainties in the fuel burn savings.
These results demonstrate that formation flight is economically beneficial in the presence of uncertainties in the fuel savings for the trailing aircraft and in the departure times of the aircraft.

In the next section, a sensitivity analysis of the solution of this experiment to the random variables that represent the departure times of the flights is carried out. The purpose of this analysis is to quantify how much uncertainty in each component of the solution is due to the different sources of uncertainty considered in the experiment.

\section{Sensitivity Analysis}
\label{sect:sensitivity_analysis}

%



In this section, a variance-based sensitivity analysis of the results obtained in Experiment B is conducted.
\textcolor{black}{The aim of the variance-based sensitivity analysis is to quantify {what proportion of} the variance of the latitude and the longitude of flight $k$ of the solution of the formation mission design problem, $\lambda_k(t, \theta)$ and  $\phi_k(t, \theta)$, respectively, is due to the variance of each departure delays of Flight 1 and Flight 2, $\theta_1$ and $\theta_2$, respectively. 
}

In this paper, this analysis relies on the computation of the so-called Sobol' indices, under the assumption that the random variables in $\theta$ are independent \cite{saltellietal:2008:gsatp}.
Sobol' indices can be directly derived from the coefficients $C_m(t)$ of the \textsf{gPC} expansion (\ref{eq:multi_dimensional_expansion_trunc}) as explained in \cite{trucchia2019}.



\begin{figure*}[ht!]
\centering
\renewcommand{\figurename}{Fig.}
\subfigure[{$S^\text{lon}_{1, 1}$ (gray), $S^\text{lon}_{1, 2}$ (black)}]{\includegraphics[width=80mm]{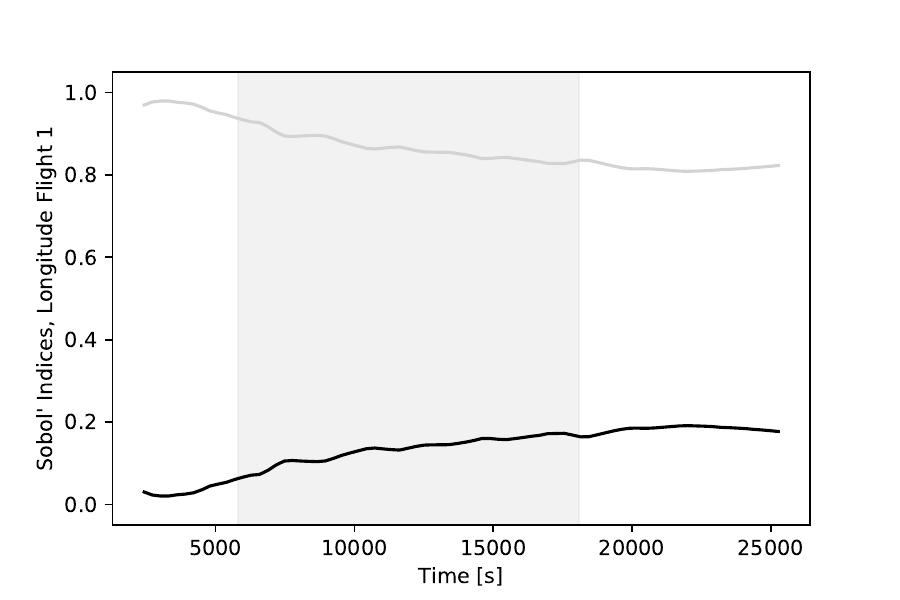}}
\subfigure[{$S^\text{lat}_{1, 1}$ (gray), $S^\text{lat}_{1, 2}$ (black)}]{\includegraphics[width=80mm]{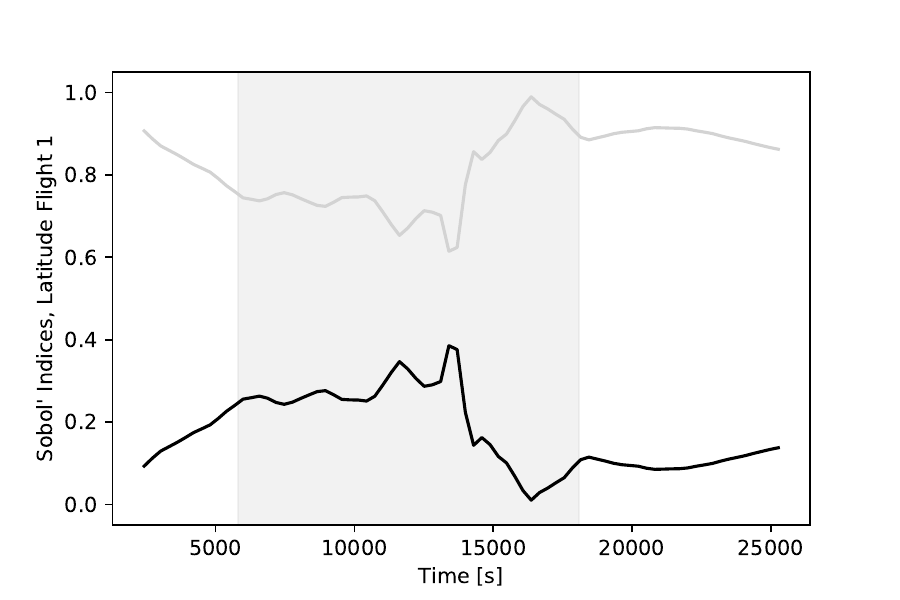}}
\subfigure[{$S^\text{lon}_{2, 1}$ (gray), $S^\text{lon}_{2, 2}$ (black)}]{\includegraphics[width=80mm]{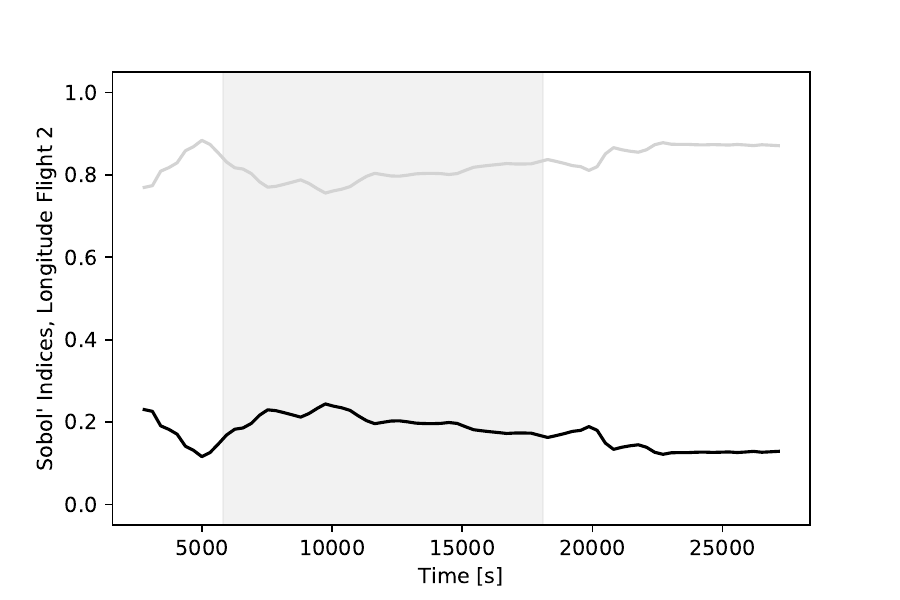}}
\subfigure[{$S^\text{lat}_{2, 1}$ (gray), $S^\text{lat}_{2, 2}$ (black)}]{\includegraphics[width=80mm]{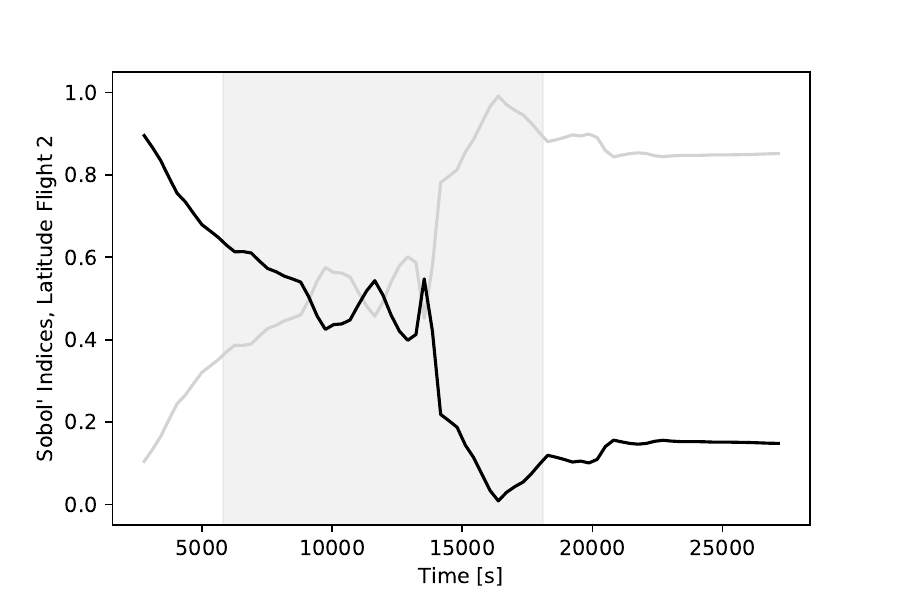}}
\caption{ {Experiment B: Sobol' indices of the geographical coordinates of the optimal trajectories of each flight. Gray and black lines correspond to the random variables $\theta_1$ and $\theta_2$, respectively.}
}
\label{fig:sobol_indices_X9_2kiwis}
\end{figure*}

The analysis of the sensitivity of the components of the solution obtained in Experiment B to the two random variables $\theta_1$ and $\theta_2$ that represent the departure delays of Flight 1 and Flight 2, respectively, is carried out employing this approach. 
The Sobol' indices of the latitude and longitude of Flight $k, k = 1,2$ with respect to the random variable $\theta_j, j=1,2$ are denoted by $S^{\text{lat}}_{k,j}$ and $S^{\text{lon}}_{k,j}$, respectively. 
They represent the proportions of the variance of the geographical coordinates of the optimal route of Flight $k$ that is due to the variance of the departure time of Flight $j$.
The Sobol' indices of the longitude and  latitude of the optimal routes obtained in the solution of Experiment B are represented as functions of time in Fig.~\ref{fig:sobol_indices_X9_2kiwis}, 
in which Sobol' indices associated with $\theta_1$ and $\theta_2$ are plotted in gray and black, respectively. 

It can be seen in Fig.~\ref{fig:sobol_indices_X9_2kiwis}.a and Fig.~\ref{fig:sobol_indices_X9_2kiwis}.b that during the whole flight, most of the variance of both geographical coordinates of the optimal route of Flight 1 is due to the variance of the random variable $\theta_1$. 
In particular, the percentage of the variance of the longitude of Flight 1 due to the variance of $\theta_1$ is around 90\%, while the percentage of the variance of the latitude of Flight 1 due to the variance of the same random variable is between 60\% and 90\% during the whole flight time.

It can be seen in Fig.~\ref{fig:sobol_indices_X9_2kiwis}.c and Fig.~\ref{fig:sobol_indices_X9_2kiwis}.d that 
during the whole flight, most of the variance of the longitude of the optimal route of Flight 2 is due to the variance of the random variable $\theta_1$. 
In contrast, the relative influence of the variance of the random variables $\theta_1$ and $\theta_2$ on the variance of the latitude of the optimal route of Flight 2 changes along the flight. 
In particular, in the first part of Flight 2, up to 8000 seconds, approximately, the variance of the random variable $\theta_2$ has the most influence on the variance of the latitude of Flight 2. 
After that, there is an intermediate part of Flight 2, between 8000 and 14000 seconds, approximately, in which the variances of $\theta_1$ and $\theta_2$ have a similar influence on the variance of the latitude of Flight 2. 
Finally, in the last part of Flight 2, the influence of the variance of $\theta_1$ on the variance of the latitude of Flight 2 becomes predominant, reaching a percentage of 90\%. 

\begin{figure}[ht!]
\centering
\renewcommand{\figurename}{Fig.}
\subfigure[{$S^\text{timing}_{1,1}$ (gray), $S^\text{timing}_{1,2}$ (black), Flight 1}]{\includegraphics[width=80mm]{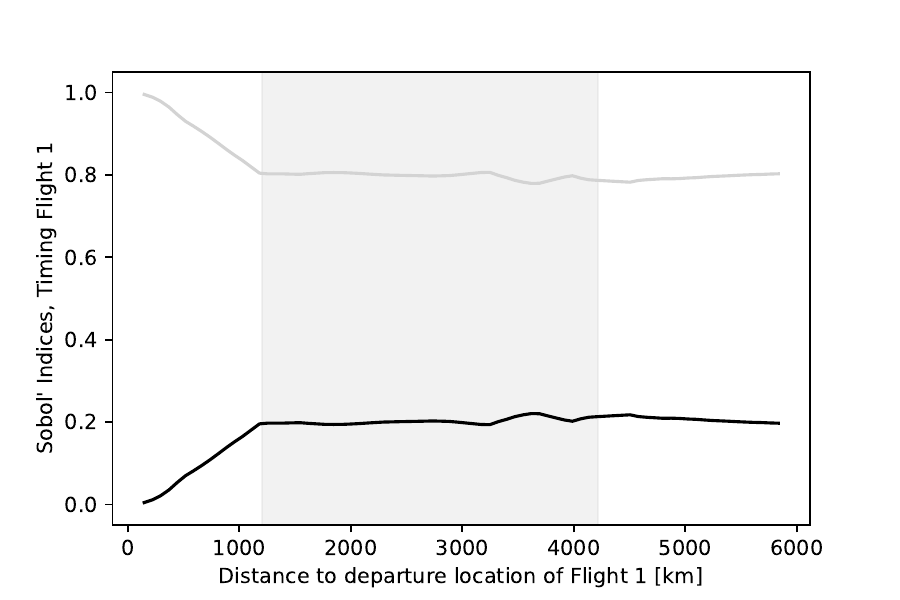}}
\subfigure[{$S^\text{timing}_{2,1}$ (gray), $S^\text{timing}_{2,2}$ (black),  Flight 2}]{\includegraphics[width=80mm]{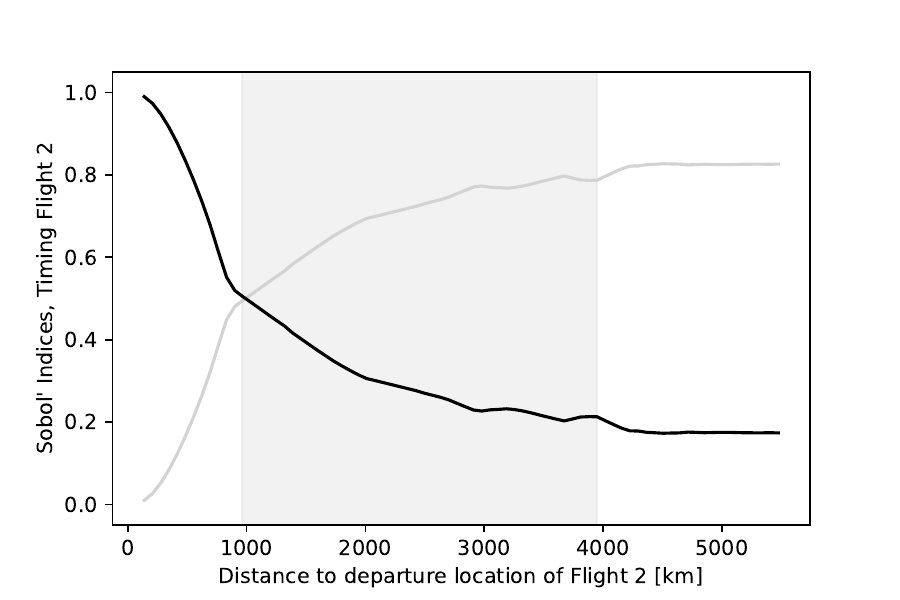}}
\caption{{Experiment B: Sobol' indices of the timing of the optimal trajectories of each flight. 
Gray and black lines correspond to the random variables $\theta_1$ and $\theta_2$, respectively.} 
}
\label{fig:sobol_indices_T_2kiwis}
\end{figure}

The Sobol' indices of the timing of Flight $k, k = 1,2$ with respect to the random variable $\theta_j, j=1,2$ are denoted by $S^\text{timing}_{k,j} (d_k)$, where $d_k$ is the orthodromic distance of flight $k$ from the departure location. They numerically quantify the proportion of the variance of the timing of the optimal route of flight $k$ that is due to the variance of $\theta_j$.
The Sobol' indices of the timing of the optimal routes obtained in the solution of Experiment B are represented as functions of the orthodromic distance in Fig.~\ref{fig:sobol_indices_T_2kiwis}, in which Sobol' indices associated with $\theta_1$ and $\theta_2$ are plotted in gray and black, respectively. 
It can be seen in this figure that, during the whole flight, most of the variance of the timing of the optimal route of Flight 1 is due to the variance of the random variable $\theta_1$, with a percentage of about 80\%. 
On the contrary, the relative influence of the variance of the random variables $\theta_1$ and $\theta_2$ on the variance of the timing of the optimal route of Flight 2 changes along the flight. 
In particular,
in the first part of Flight 2, up to the distance of 1000 km, approximately, the variance of $\theta_2$ has the most influence on the variance of the timing of the optimal route of Flight 2. After that, the influence of the variance of $\theta_1$ on the variance of the timing of the optimal route of Flight 2 becomes predominant, reaching a percentage of 80\%.

\section{Conclusions}
\label{sect:conclusions}

This paper has presented a methodology for solving the formation mission design problem for commercial aircraft in the presence of uncertainties in the fuel burn savings of the trailing aircraft and in the departure times of the aircraft. 
The solution of the problem is stochastic and the optimal aircraft trajectories are given in terms of the expected values and standard deviations of their latitude, longitude, and timing. The expected value of the direct operating costs of the formation mission is also estimated. 
Moreover, the proposed technique permits the relative contribution of each random variable to the variability of the components of the solution to be estimated.
This methodology has been applied to plan two- and three-aircraft formation missions with uncertainties.
Specifically, uncertainty in the fuel burn savings is considered in Experiment A, whereas uncertainty in the departure times of the flights is included in Experiment B.
The results of the numerical experiments indicate that uncertainties have a significant impact on the potential benefits of a formation mission. In particular, uncertainties in the departure times have a greater impact on the direct operating costs than the uncertainties in the fuel burn savings.
They reveal that the increment of the direct operating costs due to the presence of uncertainties amounts to 1.60\% and 3.76\% in Experiment A and Experiment B, respectively. 
However, even in the presence of uncertainties, expected reductions of the direct operating costs amounting to 2.16\% and 0.11\% are achieved  with formation flight compared to solo flight in Experiment A and Experiment B, respectively.
The results give interesting additional information, 
revealing, for instance, that one random variable may have predominant influence on the variability of a component of the solution. This occurs in Experiment B for the variability of the longitude of the optimal routes of both flights, the variance of which is mostly due to the variance of the random variable that represents the departure time of Flight 1. 
They also indicate that the relative influence of the random variables on the variability of a component of the solution may change during the flight. This happens in Experiment B for the variances of both the latitude and timing of Flight 2.
The results of the numerical experiments demonstrate that formation flight is economically beneficial in the presence of realistic levels of uncertainty in the fuel savings for the trailing aircraft and in the departure times of the aircraft and that the proposed methodology is an effective tool for solving the formation mission design problem for commercial aircraft in the presence of uncertainties.

\section{Acknowledgments} 
\label{sect:acknowledgments}

This work has been partially supported by the grants number
\texttt{TRA2017-91203-EXP} and \texttt{RTI2018-098471-B-C33} of the
Spanish Government.

\bibliographystyle{IEEEtran}
\bibliography{biblio_maria}

\begin{IEEEbiography}
[{\includegraphics[width=1in,height=1.5in,clip,keepaspectratio]{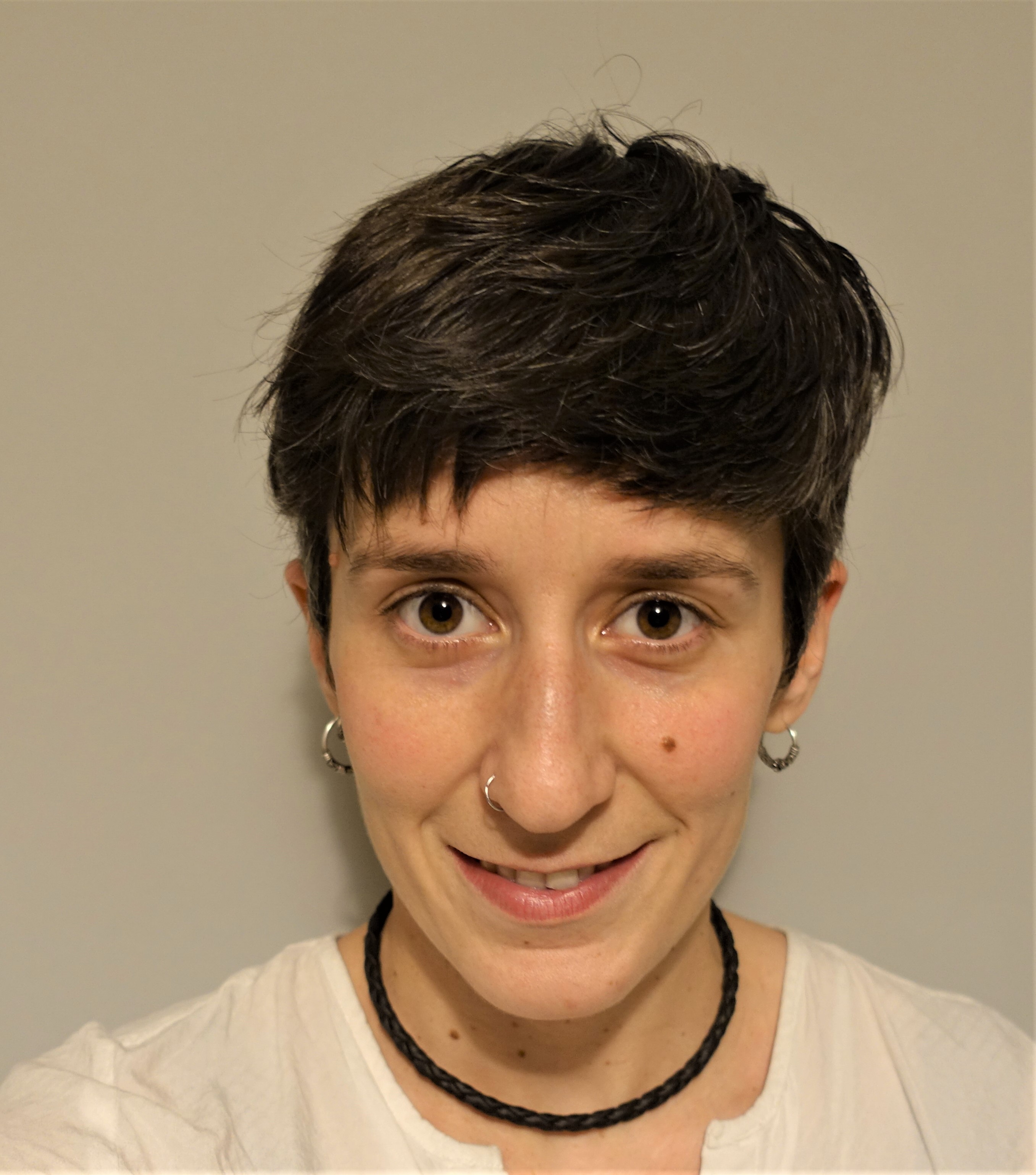}}]
{María Cerezo-Magaña} is a Teaching Assistant of {Aerospace Engineering} and a PhD student at the Universidad Rey Juan Carlos in Madrid, Spain. She received her MSc degree in Aeronautical Engineering from the Universidad Polit\'ecnica de Madrid. Her research is focused on deterministic and stochastic hybrid optimal control applied to aircraft trajectory optimization.
\end{IEEEbiography}

\begin{IEEEbiography}
[{\includegraphics[width=1in,height=1.25in,clip,keepaspectratio]{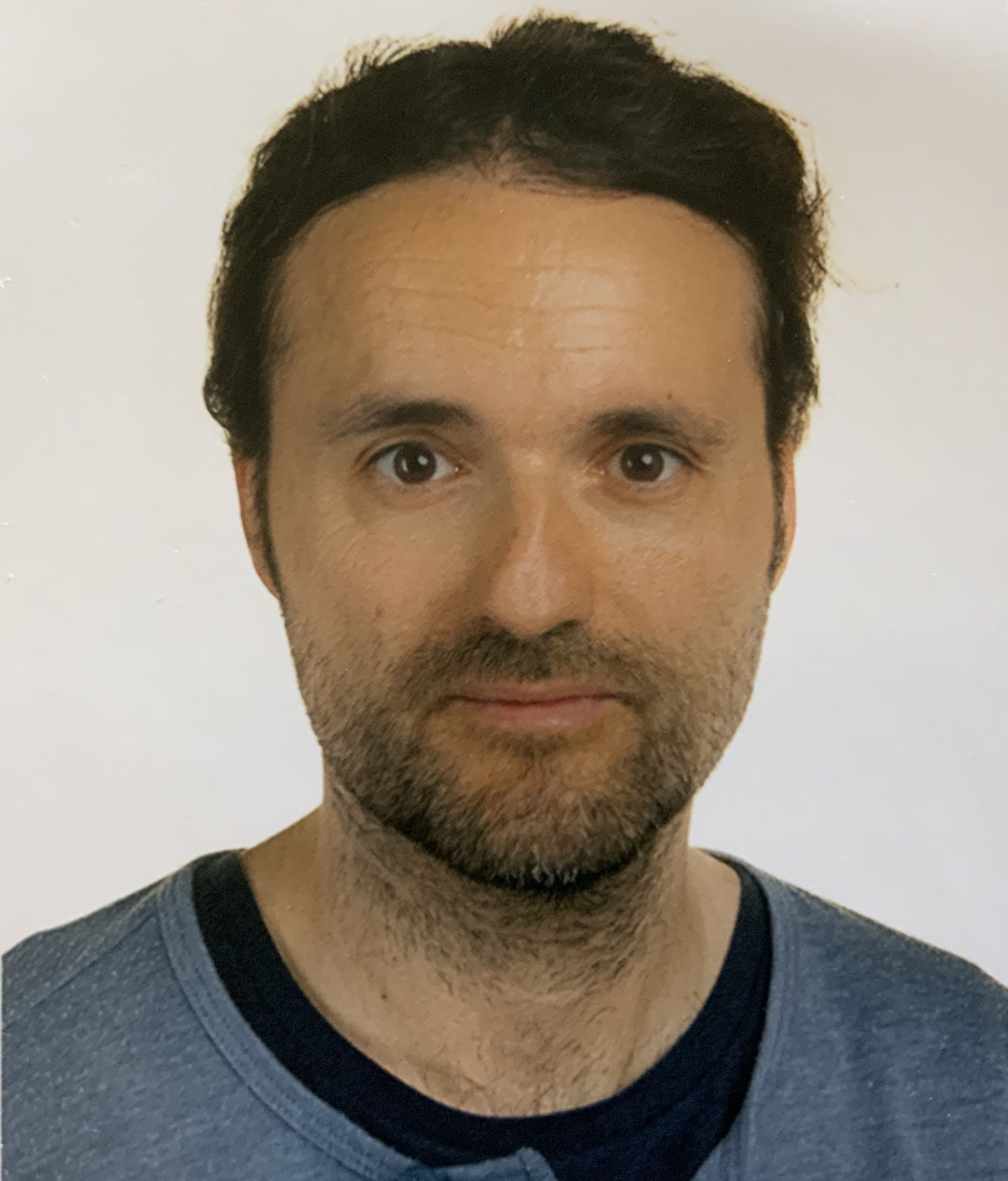}}]{Alberto Olivares}
is a Professor of Statistics and Vector Calculus at the Universidad Rey Juan Carlos in Madrid, Spain. He
received his MSc degree in Mathematics and his BSc degree in Statistics from the Universidad de Salamanca, Spain, and his PhD degree in Mathematical Engineering from the Universidad Rey Juan Carlos. He worked with the Athens University of Economics and Business. His research interests include statistical learning, stochastic hybrid optimal control and model predictive control with applications to biomedicine, robotics, aeronautics and astronautics.
\end{IEEEbiography}

\begin{IEEEbiography}
[{\includegraphics[width=1in,height=1.25in,clip,keepaspectratio]{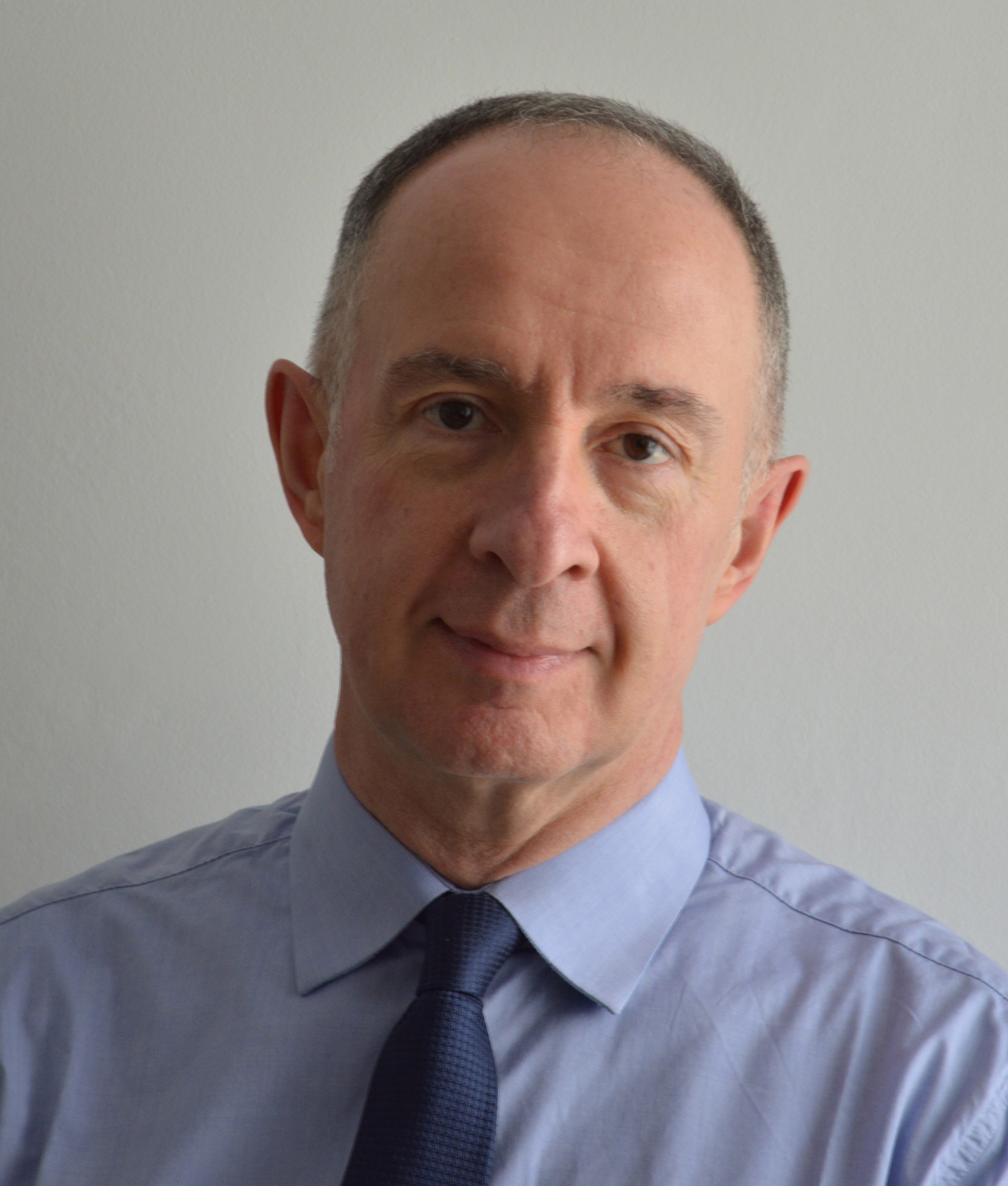}}]{Ernesto Staffetti}
is a Professor of Statistics and Control Systems at the Universidad Rey Juan Carlos in Madrid, Spain. He
received his MSc degree in Automation Engineering from the Universit\`a degli Studi di Roma ``La Sapienza,''
and his PhD degree in Advanced Automation Engineering from the Universitat Polit\`ecnica de Catalunya. 
He worked with the Universitat Polit\`ecnica de Catalunya, the
Katholieke Universiteit Leuven, the Spanish Consejo Superior de
Investigaciones Cient\'{\i}ficas, and with the University of North Carolina
at Charlotte. His research interests include stochastic hybrid optimal
control, iterative learning control and model predictive control with applications to robotics, aeronautics and astronautics.
\end{IEEEbiography}

\vfill

\end{document}